\documentclass[11pt]{article}
\usepackage{latexsym}

\font\Cp = msbm10

\newcommand{\Rrr}{\hbox{\Cp R}}
\newcommand{\Zzz}{\hbox{\Cp Z}}

\newcommand{\qed}{\mbox{$\Box$}\vspace{\baselineskip}}

\newenvironment{proof}{\noindent {\bf Proof:}}
                      {{\qed}}
\newenvironment{proof_qed}{\noindent {\bf Proof:}}
                          {\vspace{-2mm}}
\newenvironment{proof_}[1]{\noindent {\bf #1}}
                          {{\qed}}

\newtheorem{theorem}{Theorem}[section]
\newtheorem{proposition}[theorem]{Proposition}
\newtheorem{lemma}[theorem]{Lemma}
\newtheorem{definition}[theorem]{Definition}
\newtheorem{corollary}[theorem]{Corollary}

\newtheorem{examples}[theorem]{Examples}

\newtheorem{conjecture}[theorem]{Conjecture}

\newcommand{\TTch}{T}
\newcommand{\UTch}{U}
\newcommand{\FF}{C}

\newcommand{\hm}{\widehat{-1}}
\newcommand{\hz}{\widehat{0}}
\newcommand{\ho}{\widehat{1}}
\newcommand{\htwo}{\widehat{2}}

\newcommand{\tensor}{\otimes}

\newcommand{\ab}{\av\bv}
\newcommand{\av}{{\bf a}}
\newcommand{\bv}{{\bf b}}
\newcommand{\cd}{\cv\dv}
\newcommand{\ctd}{\cv\mbox{-}2\dv}
\newcommand{\cv}{{\bf c}}
\newcommand{\dv}{{\bf d}}
\newcommand{\ev}{{\bf e}}
\newcommand{\td}{{2 \dv}}
\newcommand{\xx}{{\bf x}}
\newcommand{\yy}{{\bf y}}

\newcommand{\zab}{\Zzz\langle \av,\bv \rangle}
\newcommand{\zcd}{\Zzz\langle \cv,\dv \rangle}
\newcommand{\zctd}{\Zzz\langle \cv, 2 \dv \rangle}

\newcommand{\kab}{\kk\langle \av,\bv \rangle}

\newcommand{\JH}{JH}
\newcommand{\Pyr}{{\rm Pyr}}

\newcommand{\wt}{{\rm wt}}

\newcommand{\id}{{\rm id}}

\newcommand{\kk}{{\bf k}}
\newcommand{\QSym}{{\rm QSym}}
\newcommand{\BQSym}{{\rm BQSym}}

\begin{document}

\title{The Tchebyshev Transforms of the First and Second 
Kind
}

\author{{\sc Richard EHRENBORG}
             and 
        {\sc Margaret READDY}}

\date{} 

\maketitle

\begin{abstract}
We give an
in-depth study of the Tchebyshev transforms of the first and second kind
of a poset, recently discovered by Hetyei.
The Tchebyshev transform (of the first kind)
preserves desirable combinatorial properties,
including Eulerianess (due to Hetyei) and
$EL$-shellability.
It is also a linear transformation on flag vectors.
When restricted to
Eulerian posets, it corresponds to the 
Billera, Ehrenborg and Readdy omega map $\omega$
of oriented matroids.
One consequence is that nonnegativity of the $\cd$-index is maintained
under the Tchebyshev transform.
The Tchebyshev transform of the second kind~$\UTch$
is a Hopf algebra endomorphism on 
the space of quasisymmetric functions $\QSym$.
It coincides with
Stembridge's peak enumerator $\vartheta$ 
for Eulerian posets,
but differs for general posets.
The complete
spectrum of $\UTch$ is determined,
generalizing work of Billera,
Hsiao and van Willigenburg.

\vspace{2mm}

The 
type $B$ quasisymmetric function
of a poset is introduced.
Like Ehrenborg's classical quasisymmetric function of a poset,
this map 
is a comodule morphism with
respect to the quasisymmetric functions $\QSym$.

\vspace{2mm}

Similarities among
the omega map $\omega$,
Ehrenborg's $r$-signed Birkhoff transform,
 and the Tchebyshev transforms
motivate a general study
of chain maps.
One such occurrence, 
the chain map of the second kind
$\widetilde{g}$,  is a Hopf algebra endomorphism on the quasisymmetric
functions $\QSym$
and is
an instance of Aguiar, Bergeron and Sottile's
result on
the terminal object in the category
of combinatorial Hopf algebras.
In contrast,
the chain map of the first kind
$g$ is both an algebra map 
and a comodule endomorphism on the
type $B$ quasisymmetric functions $\BQSym$.
\end{abstract}

\section{Introduction}
\setcounter{equation}{0}

The Tchebyshev transform  (of the first kind)
of a partially
ordered set, introduced by 
Hetyei~\cite{Hetyei_Tchebyshev}
and denoted by $\TTch$,
enjoys many properties.
When applied to an Eulerian poset,
this transform preserves Eulerianess~\cite{Hetyei_Tchebyshev}.
For $P$ the face lattice of a $CW$-complex,
the Tchebyshev transform yields a $CW$-complex,
that is,
the order complex of the 
Tchebyshev transform 
of a $CW$-complex triangulates the order complex of the 
original $CW$-complex~\cite{Hetyei_matrices}.
Its name derives from the fact that
when this transform is applied to the ladder poset,
the 
$\cd$-index of the resulting poset 
(expressed in terms of the variables
$\cv = \av + \bv$ and $\ev = \av - \bv$) yields
the familiar Tchebyshev polynomial of the first kind.

The $\ab$-index is a noncommutative polynomial which
encodes the flag $f$-vector of a poset.
Via a change of basis, one obtains
the $\cd$-index,
a polynomial
that removes all the linear redundancies
in the case of Eulerian posets~\cite{Bayer_Billera}.
The $\cd$-index has proven to be an extraordinarily useful
tool for studying inequalities for the face incidence
structure of 
polytopes~\cite{Billera_Ehrenborg, Ehrenborg_lifting, Ehrenborg_zonotopes}.

The omega map $\omega$, discovered by
Billera, Ehrenborg and Readdy~\cite{Billera_Ehrenborg_Readdy_om},
links the
flag $f$-vector of the intersection lattice
of a hyperplane arrangement with the
corresponding zonotope, and more generally,
the oriented matroid.
On the chain level the omega map is the inverse of a ``forgetful map''
between posets.
Aguiar and N.\ Bergeron 
observed the omega map is actually Stembridge's
peak enumerator~$\vartheta$~\cite{Stembridge}.
See~\cite{Billera_Hsiao_van_Willigenburg}
for details.

In this paper we discover new properties
of the Tchebyshev transform.  On the flag vector level
it is a linear transformation.
Surprisingly, when restricted to the class of Eulerian
posets the Tchebyshev transform is equivalent to the
omega map.
The core idea underlying this equivalence is that
the Zaslavsky's expression~\cite{Zaslavsky} for the number
of regions in a hyperplane arrangement
($\sum_{x \in P} (-1)^{\rho(x)} \cdot \mu(\hz,x)$)
applied to an Eulerian poset
gives the cardinality of the poset.
As a corollary, the Tchebyshev transform preserves 
nonnegativity of the $\cd$-index.

We also show the Tchebyshev transform preserves
$EL$-shellability.
The edge labeling we give reveals that on the flag vector level
the Tchebyshev transform of the Cartesian product of
two posets equals the dual 
diamond product of the transformed posets,
that is,
$$
\Psi(T(P \times Q)) = \Psi(T(P)) \diamond^* \Psi(T(Q)).
$$
This aforementioned proof is bijective for posets
having $R$-labelings.
A second proof is given in a more algebraic setting.
See Sections~\ref{section_Tchebyshev_Cartesian}
and~\ref{section_endomorphism}.

The theory broadens when studying
the Tchebyshev transform of the second kind $\UTch$.
(Again, Hetyei observed there
is a transform $\UTch$ which
when applied to the ladder poset
yields the Tchebyshev polynomials of the second kind.)
The Tchebyshev transform of the second kind
is a Hopf algebra endomorphism on 
the space of quasisymmetric functions $\QSym$.
The Tchebyshev transform $\UTch$ and
the peak enumerator $\vartheta$ coincide
on the $\cd$-level
but differ on the $\ab$-level,
that is, they
agree on the $\cd$-index of Eulerian posets, but differ
 on the $\ab$-index of general posets.
Billera, Hsiao and van Willigenburg~\cite{Billera_Hsiao_van_Willigenburg}
 determined the
eigenvalues and eigenvectors of Stembridge's map
$\vartheta$
when it acts on  $\zcd$ to itself.
As the transform $\UTch$ acts on $\zab$
to itself,
we can extend the diagonalization of
the map $\vartheta$ 
to this more general setting, hence deriving the complete
spectrum and eigenvectors.

There are many ways to encode the flag vector of a poset.
One is via the $\ab$-index.  Another is quasisymmetric functions.
We will now introduce a third, which we call the
{\em type $B$ quasisymmetric function} $F_{B}$ of a poset.
The type $B$ quasisymmetric functions $\BQSym$
were introduced by Chow~\cite{Chow}. 
All three encodings behave nicely under different poset products.
See Figure~\ref{figure_products}.

\begin{figure}
$$
\begin{array}{l c r c l}
\zab:       &&  \Psi(P * Q) & = & \Psi(P) \cdot \Psi(Q) \\
\QSym:      &&  F(P \times Q) & = & F(P) \cdot F(Q)     \\
\BQSym:  &&  F_{B}(P \diamond Q) & = & F_{B}(P) \cdot F_{B}(Q)     
\end{array}
$$
\caption{The product structures of
$\zab$, $\QSym$ and $\BQSym$ and their relation to poset products.}
\label{figure_products}
\end{figure}

The $\ab$-index $\Psi$ and the quasisymmetric function
$F$ of a poset are coalgebra maps.
In contrast, the type $B$ quasisymmetric function 
$F_{B}$ of a poset is a comodule map with
respect to the classical quasisymmetric function.
See Figure~\ref{figure_coproducts}.

In the study of the omega map $\omega$
relating a hyperplane arrangement
to its zonotope,
the $r$-signed Birkhoff transform $BT$ in~\cite{Ehrenborg_r-Birkhoff},
 and the Tchebyshev transforms $\TTch$ and $\UTch$,
the essential defining map has one of the following forms:
\begin{eqnarray}
 g(u) & = &
\sum_{k \geq 1}
   \sum_{u}
       \kappa(u_{(1)})
     \cdot \bv \cdot
       \widehat{g}(u_{(2)})
     \cdot \bv \cdots \bv \cdot
       \widehat{g}(u_{(k)}),       \label{equation_form_g} \\
 \widetilde{g}(u) & = &
 \sum_{k \geq 1}
   \sum_{u}
       \widehat{g}(u_{(1)})
     \cdot \bv \cdot
       \widehat{g}(u_{(2)})
     \cdot \bv \cdots \bv \cdot
       \widehat{g}(u_{(k)}) . \label{equation_form_g_tilde}
\end{eqnarray}
(The maps $\omega$, $BT$ and $\TTch$ have the form~(\ref{equation_form_g})
 and the map
$\UTch$ has the form~(\ref{equation_form_g_tilde}).
We therefore call these maps $\widetilde{g}$ and $g$
the chain maps of the first and second kind.)
This phenomenon suggests a wider theory exists on the coalgebra level.

\begin{figure}
$$
\begin{array}{l c r c l}
\zab:       &&  \Delta(\Psi(P)) & = & 
{\displaystyle
    \sum_{\hz < x < \ho} \Psi([\hz,x]) \tensor \Psi([x,\ho]) } \\[6 mm]
\QSym:      &&  \Delta^{\QSym}(F(P)) & = &
{\displaystyle
    \sum_{\hz \leq x \leq \ho} F([\hz,x]) \tensor F([x,\ho]) } \\[6 mm]
\BQSym:  &&  \Delta^{\BQSym}(F_{B}(P)) & = &
{\displaystyle
    \sum_{\hz < x \leq \ho} F_{B}([\hz,x]) \tensor F([x,\ho]) }
\end{array}
$$
\caption{The coalgebra structures of
$\zab$, $\QSym$ and $\BQSym$ and their relation to posets.}
\label{figure_coproducts}
\end{figure}

In Sections~\ref{section_chain_maps}
and~\ref{section_quasisymmetric_type_B}
we study general functions of these types.
We show
the chain map of the second kind
$\widetilde{g}$ is a Hopf algebra endomorphism on quasisymmetric
functions.
This is a concrete example of
Aguiar, Bergeron and Sottile's
theorem 
that the algebra of quasisymmetric functions $\QSym$
is the terminal object in the category
of combinatorial
Hopf algebras~\cite{Aguiar_Bergeron_Sottile}.
Furthermore,
the chain map of the first kind
$g$ is an algebra map on the type $B$ quasisymmetric
functions.
See Theorems~\ref{theorem_widetilde_M} and~\ref{theorem_algebra}.
The map $g$ is also a comodule endomorphism on 
the type $B$ quasisymmetric functions.
See Theorem~\ref{theorem_g_comodule_endomorphism}.

We end the paper with concluding remarks and
many questions for further study.

\section{Background Definitions}
\setcounter{equation}{0}

For a graded poset $P$ with minimal element $\hz$ and maximal
element $\ho$, let
$P \cup \{ \hm \}$ and
$P \cup \{ \htwo \}$ denote
$P$ adjoined with a new minimal
element $\hm$, respectively
a new
maximal element $\htwo$.
For a chain
$c = \{\hz = x_0 < x_1 < \cdots < x_k = \ho\}$
in $P$ 
define the {\em weight} of the chain $c$
by
$$    \wt(c)
    =
         (\av-\bv)^{\rho(x_{0},x_{1}) - 1}
        \cdot
         \bv
        \cdot
         (\av-\bv)^{\rho(x_{1},x_{2}) - 1}
        \cdot
         \bv
        \cdots
         \bv
        \cdot
         (\av-\bv)^{\rho(x_{k-1},x_{k}) - 1} , $$
where 
$\av$ and $\bv$ are noncommutative variables.
The
{\em \ab-index}
of the poset $P$ is defined as
$$
	\Psi(P) = \sum_{c} \wt(c),
$$
where the sum is over all chains $c$ in $P$.

A poset is {\em Eulerian} if every interval $[x,y]$, where $x < y$,
has the same number of elements of even rank as elements of odd rank.
When $P$ is Eulerian
the $\ab$-index of $P$
can be written in terms of 
$\cv = \av + \bv$ and
$\dv = \av \cdot \bv + \bv \cdot \av$;
see~\cite{Bayer_Klapper}.
The resulting noncommutative polynomial 
is called the 
{\em $\cd$-index}.
Its importance lies in that it removes all the linear redundancies
in the flag $f$-vector entries~\cite{Bayer_Billera},
geometric operations on a polytope translate as
operators of the corresponding 
$\cd$-index~\cite{Ehrenborg_Johnston_Rajagopalan_Readdy,
Ehrenborg_Readdy_c},
and
is amenable to algebraic
techniques to derive
inequalities
on the flag
vectors~\cite{Billera_Ehrenborg,Ehrenborg_lifting,Ehrenborg_zonotopes}.

On the ring $\zab$ define a coproduct $\Delta$
by defining it on an $\ab$-monomial
$u_{1} u_{2} \cdots u_{k}$ by
$$
     \Delta(u_{1} u_{2} \cdots u_{k})
   =
     \sum_{i=1}^{k} u_{1} \cdots u_{i-1} \tensor u_{i+1} \cdots u_{k} ,
$$
where each $u_{i}$ is either an $\av$ or a $\bv$.
It is straightforward to verify the coproduct $\Delta$
is coassociative, that is,
$(\id \tensor \Delta) \circ \Delta
    =
 (\Delta \tensor \id) \circ \Delta$.
Hence define
$\Delta^{k} : \zab \longrightarrow \zab^{\tensor k}$
by
$\Delta^{1} = \id$
and
$\Delta^{k+1} = (\id \tensor \Delta^{k}) \circ \Delta$.
The coproduct $\Delta$ satisfies the
Newtonian condition:
\begin{equation}
        \Delta(u \cdot v)
     =
        \sum_{u} u_{(1)} \tensor u_{(1)} v
     +
        \sum_{v} u v_{(1)} \tensor v_{(1)}     .
\label{equation_Newton}
\end{equation}

The essential property of the coproduct $\Delta$
is that
it makes the $\ab$-index into a coalgebra
homomorphism~\cite{Ehrenborg_Readdy_c}.
\begin{theorem}
For a graded poset $P$ we have
$$
     \Delta(\Psi(P))
   =
     \sum_{\hz < x < \ho} \Psi([\hz,x]) \tensor \Psi([x,\ho])  .
$$
\end{theorem}
This allows for computations 
on posets to be translated into
the coalgebra $\zab$;
see~\cite{Billera_Ehrenborg,Billera_Ehrenborg_Readdy_om,Ehrenborg_r-Birkhoff,Ehrenborg_Fox,Ehrenborg_Johnston_Rajagopalan_Readdy,Ehrenborg_Readdy_c}.

Using the coassociativity, we have the following
corollary.
\begin{corollary}
\begin{equation}
     \Delta^{k}(\Psi(P))
   =
     \sum_{\hz = x_{0} < x_{1} < \cdots < x_{k} = \ho}
          \Psi([x_{0},x_{1}])
        \tensor
          \Psi([x_{1},x_{2}])
        \tensor
          \cdots
        \tensor
          \Psi([x_{k-1},x_{k}])   .
\label{equation_k_chains}
\end{equation}
\end{corollary}

There is an involution on $\zab$ that sends each
monomial $u = u_{1} u_{2} \cdots u_{k}$ to its
reverse $u^{*} = u_{k} \cdots u_{2} u_{1}$.
Directly we have
$(u \cdot v)^{*} = v^{*} \cdot u^{*}$,
$\Delta(u^{*}) = \sum_{u} u_{(2)}^{*} \tensor u_{(1)}^{*}$
and
$\Psi(P^{*}) = \Psi(P)^{*}$
where
$P^{*}$
denotes the dual of the poset $P$.

\section{Quasisymmetric functions}
\label{section_quasisymmetric}
\setcounter{equation}{0}

Another way to encode the flag $f$-vector of a poset $P$
is by the quasisymmetric function~$F(P)$;
see~\cite{Ehrenborg_Hopf}.
Let $P$ be a poset of rank $n$, where $n \geq 0$.
The {\em quasisymmetric function} of the poset $P$ is
defined as the limit
$$   F(P)
  =
     \lim_{m \longrightarrow \infty}
        \sum_{\hz = x_{0} \leq x_{1} \leq \cdots \leq x_{m} = \ho}
             t_{1}^{\rho(x_{0},x_{1})}
           \cdot
             t_{2}^{\rho(x_{1},x_{2})}
           \cdots
             t_{m}^{\rho(x_{m-1},x_{m})}  .  $$
Observe that for $m = 2$ this sum is a homogeneous
rank-generating function, that is, it encodes the $f$-vector
of the poset. For larger $m$ it encodes all the entries
in the flag $f$-vector of cardinality less than or equal to $m-1$.

The polynomial $F(P)$ is homogeneous of degree $n$ in the infinitely-many
variables $t_{1}, t_{2}, \ldots$.  It also enjoys the following
quasisymmetry: for $i_{1} < i_{2} < \cdots < i_{k}$
and $j_{1} < j_{2} < \cdots < j_{k}$
the coefficients of
$t_{i_{1}}^{p_{1}} \cdot t_{i_{2}}^{p_{2}} \cdots t_{i_{k}}^{p_{k}}$
and
$t_{j_{1}}^{p_{1}} \cdot t_{j_{2}}^{p_{2}} \cdots t_{j_{k}}^{p_{k}}$
are the same. Polynomials in the variables $t_{1}, t_{2}, \ldots$
are called {\em quasisymmetric} and the algebra of these
polynomials are denote by $\QSym$. It is straightforward
to observe that a linear basis for $\QSym$ is given
by the {\em monomial quasisymmetric function},
defined by
$$        M_{(p_{1}, p_{2}, \ldots, p_{k})}
  =
    \sum_{i_{1} < i_{2} < \cdots < i_{k}}
         t_{i_{1}}^{p_{1}}
         t_{i_{2}}^{p_{2}}
      \cdots
         t_{i_{k}}^{p_{k}}     .  $$

Define a linear map $\gamma$ from $\kab$ to $\QSym$
by
$$ \gamma\left(
        (\av-\bv)^{p_{1}-1}
      \cdot
        \bv
      \cdot
        (\av-\bv)^{p_{2}-1}
      \cdot
        \bv
      \cdots
        \bv
      \cdot
        (\av-\bv)^{p_{k}-1} \right)
 =
        M_{(p_{1}, \ldots, p_{k})}     . $$
This map is an isomorphism between $\kab$ and
quasisymmetric functions having no constant term.
For a poset $P$ of rank greater than or equal to
one, we have  $\gamma(\Psi(P)) = F(P)$.
For the one element poset $\bullet$ of rank $0$,
let $F(\bullet) = 1_{\QSym}$.
Here we write $1_{\QSym}$ for the identity
element of the quasisymmetric functions
in order to distinguish it from the unit
in $\zab$.
For more on the Hopf algebra structure of the quasisymmetric
functions~$\QSym$, we refer the reader to~\cite{Ehrenborg_Hopf}.

Let us mention two important identities for
the quasisymmetric function $F(P)$ of a graded poset~$P$.
For $P$ and $Q$ two graded posets, we have
\begin{eqnarray}
F(P) \cdot F(Q)
  & = &
F(P \times Q) , 
\label{equation_product} \\
\Delta^{\QSym}(F(P))
  & = &
 \sum_{\hz \leq x \leq \ho}
  F([\hz,x]) \tensor F([x,\ho])  \nonumber \\
  & = &
 F(P) \tensor 1_{\QSym} + 1_{\QSym} \tensor F(P) +
 \sum_{\hz < x < \ho}
  F([\hz,x]) \tensor F([x,\ho])  ,
\label{equation_difference}
\end{eqnarray}
where equation~(\ref{equation_difference})
is valid when
the poset $P$ has rank at least $1$.
Note the coproduct on quasisymmetric functions
differs from the coproduct on $\ab$-polynomials.
In order to avoid confusion, we are denoting the coproduct
on quasisymmetric functions by
$\Delta^{\QSym}(f) = \sum^{\QSym}_{f} f_{(1)} \tensor f_{(2)}$.
For proofs of these identities,
see~\cite[Proposition~4.4]{Ehrenborg_Hopf}.

{}From a poset perspective identities~(\ref{equation_product})
and (\ref{equation_difference})
define the algebra and coalgebra structure of
the quasisymmetric functions~$\QSym$.
Equation~(\ref{equation_difference})
also motivates the following relation between the two
coproducts~$\Delta$ and~$\Delta^{\QSym}$:
\begin{equation}
     \Delta^{\QSym}(\gamma(v))
        =
     \gamma(v) \tensor 1_{\QSym}
        +
     1_{\QSym} \tensor \gamma(v)
        +
     \sum_{v} \gamma(v_{(1)}) \tensor \gamma(v_{(2)}) .  
\end{equation}

\section{Enumerating flags in the Tchebyshev transform of a poset}
\setcounter{equation}{0}

\begin{definition}
For a graded poset $P$ define the 
Tchebyshev transform (of the first kind) $T(P)$ to be the graded poset 
with elements given by the set 
$$   T(P)  
  =
     \{ [x,y] \:\: : \:\: x,y \in P \cup \{\hm\}, \:\: x < y  \}, $$
and the cover relation given by 
the following three rules:
\begin{enumerate}
\item[(i)] $[x,y] \prec_{T(P)} [y,w]$ if $y \prec w$,

\item[(ii)] $[x,y] \prec_{T(P)} [x,w]$ if $y \prec w$, and

\item[(iii)] $[x,\ho] \prec_{T(P)}  \ho_{T(P)}$.
\end{enumerate}
\end{definition}
As a remark,
Hetyei's original definition of the Tchebyshev transform is
in terms of the order relation rather than the cover relation
of the poset.
Note the  rank function $\rho_{T(P)}$ on $T(P)$
satisfies
$\rho_{T(P)}([x,y]) = \rho_{P}(y)$
and
$\rho_{T(P)}(\ho_{T(P)}) = \rho(P) + 1$.

Our interest in studying the Tchebyshev transform of posets
arises from the following 
surprising result of Hetyei~\cite{Hetyei_Tchebyshev}.
\begin{theorem}
Let $P$ be an Eulerian poset.
Then the Tchebyshev transform of $P$ is also an Eulerian poset.
\label{theorem_Tchebyshev_Eulerian}
\end{theorem}

We now prove a proposition which can been viewed as an analogue of
a result of Bayer and Sturmfels~\cite{Bayer_Sturmfels}
(see Proposition~4.6.2 in~\cite{BLSWZ})
and
of Proposition~4.1 in~\cite{Ehrenborg_r-Birkhoff}.
This connection will be made clearer in
Sections~\ref{section_chain_maps}
and~\ref{section_quasisymmetric_type_B}.

\begin{proposition}
For a chain
$c = \{\hz = x_0 < x_1 < \cdots < x_k = \htwo\}$ in $P \cup \{\htwo\}$,
the cardinality of the inverse image of $c$ is given by
$$    |z^{-1}(c)|
    =
      \prod_{i=1}^{k-1}
            | [x_{i-1}, x_{i}] |
   .  $$
\label{proposition_inverse_image}
\end{proposition}

To prove this proposition we need a lemma
and its corollary.

\begin{lemma}
Given three elements $x < y < w$ in the poset $P$,
the condition
 $[x,y] <_{T(P)} [z,w]$ is equivalent to
either
$z = x$ 
or
$z \in [y,w)$. 
\end{lemma}
\begin{proof}
We proceed by induction on $\rho(y,w)$.
If
$\rho(y,w) = 1$, we have by definition 
that the element $z$ is either $x$
(condition ($ii$))
or
$y$ (condition ($i$)).
Assume now that
$\rho(y,w) \geq 2$ and
let
$[u,v]$ be an atom in the interval
$[[x,y],[z,w]]$.
Since
$\rho(v,w) < \rho(y,w)$ we have by the
induction hypothesis that either
$z = u$ or $z \in [v,w)$.  The union of all such intervals
$[v,w)$ is the open interval $(y,w)$.
Moreover, since $v$ covers~$y$, we have that
$u$ is either $x$ or $y$.
That is, the only choices for $z$ are
$\{x\} \cup \{y\} \cup (y,w) = \{x\} \cup [y,w)$,
proving the induction step.
\end{proof}

\begin{corollary}
Given three elements $x < y < w$ in the poset $P$,
the number of elements $z$ such that
$[x,y] <_{T(P)} [z,w]$ equals
the cardinality of the interval $[y,w]$.
\label{corollary_x_y_w}
\end{corollary}

For a graded poset $P$
let $z: T(P) \rightarrow P \cup \{\htwo\}$
be the 
map
$z([x,y]) = y$ and
$z(\ho_{T(P)}) = \htwo$.
Observe the map $z$ is order and rank preserving 
and hence preserves chains and the
weight of chains.

The proof of Proposition~\ref{proposition_inverse_image}
follows by repeated use of Corollary~\ref{corollary_x_y_w}.

\section{The Tchebyshev transform on $\ab$-polynomials}
\label{section_Tchebyshev_transform}
\setcounter{equation}{0}

In this section we express the $\ab$-index of the Tchebyshev
transform in terms of the $\ab$-index of the original
poset.

Define two linear maps $A$ and $\FF$ from $\zab$ to $\Zzz$
as follows.
The map $A$ is the algebra map with
$A(\av) = 1$ and $A(\bv) = 0$ and the map $\FF$ is given by
the relation
$$   \FF(u) = 2 \cdot A(u) + \sum_{u} A(u_{(1)}) \cdot A(u_{(2)}) .  $$
Here we are using the usual Sweedler notation~\cite{Sweedler}.
\begin{lemma}
For a graded poset $P$ we have 
$A(\Psi(P)) = 1$ and
$\FF(\Psi(P))$ is the cardinality of the poset~$P$.
\end{lemma}
\begin{proof}
The first identity
$A(\Psi(P)) = 1$
was already observed in~\cite{Billera_Ehrenborg_Readdy_om}.
The second identity follows from
\begin{eqnarray*}
\FF(\Psi(P))
  & = &
     2 \cdot A(\Psi(P))
   +
     \sum_{\hz < x < \ho } A(\Psi([\hz,x])) \cdot A(\Psi([x,\ho])) \\
  & = &
     2 
   +
     \sum_{\hz < x < \ho } 1    \\
  & = &
   | P |    ,  
\end{eqnarray*}
where the first step follows from the fact the $\ab$-index is
a coalgebra homomorphism.
\end{proof}

\begin{lemma}
The linear map $\FF$ satisfies the recursion
\begin{eqnarray*}
\FF(1) & = & 2 , \\
\FF(\av \cdot u) & = & A(u) + \FF(u) , \\
\FF(\bv \cdot u) & = & A(u) . \\
\end{eqnarray*}
\end{lemma}
\begin{proof_qed}
Directly
$\FF(1)  =  2 \cdot A(1)  =  2$.
For the second identity we have
by the Newtonian condition~(\ref{equation_Newton})
\begin{eqnarray*}
\FF(\av \cdot u)
  & = &
2 \cdot A(\av \cdot u) 
  +
A(1) \cdot A(u) 
  +
\sum_{u} A(\av \cdot u_{(1)}) \cdot A(u_{(2)})  \\
  & = &
3 \cdot A(u) 
  +
\sum_{u} A(u_{(1)}) \cdot A(u_{(2)})  \\
  & = &
A(u) 
  +
\FF(u) . 
\end{eqnarray*}
Similarly, the third identity follows from
\begin{eqnarray*}
\hspace*{25 mm}
\FF(\bv \cdot u)
  & = &
2 \cdot A(\bv \cdot u) 
  +
A(1) \cdot A(u) 
  +
\sum_{u} A(\bv \cdot u_{(1)}) \cdot A(u_{(2)})  \\
  & = &
A(u)    .
\hspace*{75 mm}
\hspace*{25 mm}
\qed
\end{eqnarray*}
\end{proof_qed}

We now consider three linear operators on $\zab$.
For a homogeneous $\ab$-polynomial $u$
define $\kappa$ and $\nu$ by
$$    \kappa(u) = A(u) \cdot (\av - \bv)^{\deg(u)}
    \:\:\:\: \mbox{ and } \:\:\:\:
      \nu(u) = \FF(u) \cdot (\av - \bv)^{\deg(u)}   ,  $$
and extend by linearity.
Define $\TTch$ by the sum
\begin{equation}
  \TTch(u) = \sum_{k \geq 1}
             \sum_{u}
             \nu(u_{(1)}) \cdot \bv \cdot
             \nu(u_{(2)}) \cdot \bv \cdots
                               \bv \cdot
             \nu(u_{(k-1)}) \cdot \bv \cdot
             \kappa(u_{(k)})   ,
\label{equation_TTch}
\end{equation}
where the coproduct is into $k$ parts.

The slight abuse of notation between the
Tchebyshev transform of a graded poset
and the Tchebyshev transform of $\ab$-monomials
is explained by the following theorem.

\begin{theorem}
The $\ab$-index of the Tchebyshev transform of 
a graded poset $P$
is given by
$$     \Psi(T(P)) = \TTch(\Psi(P) \cdot \av)   .  $$
\end{theorem}
\begin{proof}
Using the chain definition of the $\ab$-index 
and Proposition~\ref{proposition_inverse_image},
we have
\begin{eqnarray*}
\Psi(T(P))
  & = &
\sum_{c} |z^{-1}(c)| \cdot \wt(c) \\
  & = &
  \sum_{k \geq 1}
  \sum_{c} 
      \prod_{i=1}^{k-1}
         \FF(\Psi([x_{i-1},x_{i}]))  
        \cdot
      A(\Psi([x_{k-1},x_{k}]))
    \cdot \wt(c)    \\
  & = &
  \sum_{k \geq 1}
  \sum_{c} 
      \prod_{i=1}^{k-1}
         \FF(\Psi([x_{i-1},x_{i}]))  
        \cdot
      A(\Psi([x_{k-1},x_{k}])) \\
  &   &
    \cdot 
\left(
      \prod_{i=1}^{k-1}
         (\av-\bv)^{\rho(x_{i-1},x_{i}) - 1}
        \cdot
         \bv
\right)
        \cdot
         (\av-\bv)^{\rho(x_{k-1},x_{k}) - 1} \\
    & = &
  \sum_{k \geq 1}
  \sum_{c} 
\left(
      \prod_{i=1}^{k-1}
         \nu([x_{i-1},x_{i}])
        \cdot
         \bv
\right)
        \cdot
         \kappa([x_{k-1},x_{k}])   \\
    & = &
  \sum_{k \geq 1}
  \sum_{u} 
\left(
      \prod_{i=1}^{k-1}
         \nu(u_{(i)})
        \cdot
         \bv
\right)
        \cdot
         \kappa(u_{(k)})   \\
    & = &
 \TTch(u) .
\end{eqnarray*}
Here the second to last step uses the fact the $\ab$-index
is a coalgebra homomorphism
and $u$ is the
$\ab$-index of $P \cup {\htwo}$, that is,
$u = \Psi(P \cup {\htwo}) = \Psi(P) \cdot \av$.
\end{proof}

\begin{proposition}
The operator $\TTch$ satisfies the following functional
identity:
$$    \TTch(u)  =   \kappa(u)
                 +
                   \sum_{u}  \nu(u_{(1)}) \cdot \bv \cdot \TTch(u_{(2)}) . $$
\label{proposition_functional_equation}
\end{proposition}

\section{Connection with the $\omega$ operator of oriented matroids}
\label{section_omega}
\setcounter{equation}{0}

We begin by recalling the $\omega$ map for oriented 
matroids~\cite{Billera_Ehrenborg_Readdy_om}.

\begin{theorem}
Let $\omega: \zab \rightarrow \zctd$ be the linear map 
defined on
monomials
in the variables
$\av$ and
$\bv$
by replacing each occurrence of
$\ab$ by $\td$
and the remaining letters
with
$\cv$'s.
Let $R$ be the lattice of regions 
and $L$ be the lattice of flats 
of an oriented matroid.
Then the $\cd$-index of 
$R$ is given by
$$
     \Psi(R) = \omega(\av \cdot \Psi(L))^*.
$$
In fact, the $\cd$-index of the lattice of regions $R$
is indeed a $\ctd$-index.
\end{theorem}

Hsiao has found an analogous version of this theorem for the
Birkhoff transform of a distributive lattice~\cite{Hsiao}.
Ehrenborg has generalized Hsiao's work to an $r$-signed
Birkhoff transform~\cite{Ehrenborg_r-Birkhoff}.
In this section we show the Tchebyshev transform is
likewise connected to the omega map.
This allows us to conclude the Tchebyshev transform preserves
nonnegativity of the $\cd$-index.

\begin{theorem}
For $\cd$-polynomials $v$ we have
$$      \TTch(v \cdot \av)
    =
        \omega(\av \cdot v^{*})^{*}   .  $$
\label{theorem_TTch_omega}
\end{theorem}
\begin{proof}
Following~\cite{Billera_Ehrenborg_Readdy_om}
let $\eta$ be the unique operator on $\zab$ such that
$$ \eta(\Psi(P))
 =
   \left(
      \sum_{\hz \leq x \leq \ho}
             (-1)^{\rho(x)}
                \cdot
             \mu(\hz,x)
   \right)
      \cdot
   (\av-\bv)^{\rho(P)-1}    ,  $$
for all posets $P$.
Next, let the operator $\varphi$ be defined as follows:
\begin{equation}
   \varphi(u) = \sum_{k \geq 1}
             \sum_{u}
             \kappa(u_{(1)}) \cdot \bv \cdot
             \eta(u_{(2)}) \cdot \bv \cdots
                               \bv \cdot
             \eta(u_{(k)})                      .
\label{equation_phi}
\end{equation}
By Proposition~5.5
in~\cite{Billera_Ehrenborg_Readdy_om}
we have 
$\omega(\av \cdot v) = \varphi(\av \cdot v)$.
Also observe 
\begin{equation}
  \TTch(u^{*})^{*}
           =
             \sum_{k \geq 1}
             \sum_{u}
             \kappa(u_{(1)}) \cdot \bv \cdot
             \nu(u_{(2)}) \cdot \bv \cdots
                               \bv \cdot
             \nu(u_{(k)})   .
\label{equation_TTch_star}
\end{equation}

For any Eulerian poset $P$ we have 
$\eta(\Psi(P)) = \nu(\Psi(P))$
since 
$(-1)^{\rho(x)} \cdot \mu(\hz,x) = 1$ for all
elements $x$ in an Eulerian poset $P$.
Since the $\cd$-indexes of all Eulerian posets span all
$\cd$-polynomials, we have for all $\cd$-polynomials
$v$ that
$   \eta(v) = \nu(v)$.
Now consider the coproduct $\Delta^{k}$ applied
to $u = \av \cdot v$, where $v$ is a $\cd$-polynomial.
We obtain
$$  \Delta^{k}(\av \cdot v) 
        \in
    \zab \tensor \zcd^{\tensor(k-1)}   .  $$
Hence the expressions in 
equations~(\ref{equation_phi})
and~(\ref{equation_TTch_star})
agree on $u = \av \cdot v$.
\end{proof}

Recall that Hetyei proved the Tchebyshev transform
preserves Eulerianness.
See Theorem~\ref{theorem_Tchebyshev_Eulerian}.
We obtain two important corollaries.
\begin{theorem}
If an Eulerian poset $P$ has a non-negative $\cd$-index
so does the Tchebyshev transform~$T(P)$,
that is, $\Psi(P) \geq 0$ implies $\Psi(T(P)) \geq 0$.
\label{theorem_non-negativity}
\end{theorem}
\begin{proof}
The $\cd$-polynomial $\Psi(P)$ has non-negative terms
as an $\ab$-polynomial.
Applying Theorem~\ref{theorem_TTch_omega}
and observing
that $\omega$ sends an $\ab$-monomial
to a $\ctd$-monomial, we see that non-negativity
is preserved.
\end{proof}

\begin{corollary}
The Tchebyshev transform $T(P)$ of an Eulerian poset $P$
has a $\ctd$-index, that is, the $\cd$-index $\Psi(T(P))$ belongs
to $\zctd$.
\label{corollary_ctd}
\end{corollary}

Since a given $\cd$-monomial expands into $2^{k}$ $\ab$-monomials,
where $k$ is the number of $\cv$'s and $\dv$'s
appearing in the monomial, we also have:
\begin{corollary}
Let $u$ be a $\cd$-monomial consisting of $k$ letters.
Then the Tchebyshev transform $\TTch(u \cdot \av)$
is a sum of $2^{k}$ $\ctd$-monomials.
\label{corollary_2_power}
\end{corollary}

Recall the hyperplane arrangement in $\Rrr^n$ consisting of the
$n$ coordinate hyperplanes 
$x_i = 0$ for $1 \leq i \leq n$ 
has intersection lattice corresponding to the Boolean algebra
$B_n$.
The regions of this arrangement 
correspond to
the $n$-dimensional crosspolytope $C_n$.
Hence another corollary of Theorem~\ref{theorem_TTch_omega} is:
\begin{corollary}
\label{corollary_Tchebyshev_the_Boolean_algebra}
The $\cd$-index of the Tchebyshev transform of the Boolean algebra $B_{n}$
is given by the $\cd$-index
of the $n$-dimensional crosspolytope $C_{n}$,
that is,
$$
     \Psi(T(B_{n})) = \Psi(C_{n}).
$$
\end{corollary}

\section{Recursions for the Tchebyshev transform}
\setcounter{equation}{0}

In this section we develop recursions for computing
the Tchebyshev transform. They are especially important
for $\ab$-polynomials.

Define a new operator 
$\sigma$ on $\zab$
by
$$   \sigma(u)
   =
     \sum_{u} \kappa(u_{(1)}) \cdot \bv \cdot \TTch(u_{(2)})   ,  $$
where $u$ is an $\ab$-polynomial.

\begin{proposition}
The operator $\TTch$ satisfies the following joint recursion
with the operator $\sigma$:
\begin{eqnarray}
\TTch(1)
  & = &
1 \\
\TTch(\av \cdot u)
  & = & 
(\av + \bv) \cdot \TTch(u)
   +
(\av - \bv) \cdot \sigma(u) , \label{equation_tau_a} \\
\TTch(\bv \cdot u)
  & = & 
2 \bv \cdot \TTch(u) 
   +
(\av-\bv) \cdot \sigma(u)   , \label{equation_tau_b} \\
\sigma(1)
  & = &
0 ,  \\
\sigma(\av \cdot u)
  & = & 
\bv \cdot \TTch(u)
    +
(\av-\bv) \cdot \sigma(u) ,  \\
\sigma(\bv \cdot u)
  & = &
\bv \cdot \TTch(u) . 
\end{eqnarray}
\label{proposition_ab_recursion}
\end{proposition}
\begin{proof_qed}
Directly $\TTch(1) = 1$ and $\sigma(1) = 0$.
Using the Newtonian condition~(\ref{equation_Newton}),
we have
\begin{eqnarray*}
\TTch(\av \cdot u)
  & = &
\kappa(\av \cdot u)
   +
\nu(1) \cdot \bv \cdot \TTch(u) 
   +
\sum_{u}  \nu(\av \cdot u_{(1)}) \cdot \bv \cdot \TTch(u_{(2)}) \\
  & = &
(\av - \bv) \cdot \kappa(u)
   +
2 \bv \cdot \TTch(u) 
   +
(\av-\bv) \cdot \sum_{u}
      (\kappa(u_{(1)}) + \nu(u_{(1)})) \cdot \bv \cdot \TTch(u_{(2)}) \\
  & = &
(\av + \bv) \cdot \TTch(u)
   +
(\av-\bv) \cdot \sum_{u} \kappa(u_{(1)}) \cdot \bv \cdot \TTch(u_{(2)}) \\
  & = &
(\av + \bv) \cdot \TTch(u)
   +
(\av - \bv) \cdot \sigma(u) .
\end{eqnarray*}
Here we have used 
the functional equation in Proposition~\ref{proposition_functional_equation}
in the first and third equalities.
Similarly, we obtain
\begin{eqnarray*}
\TTch(\bv \cdot u)
  & = &
\kappa(\bv \cdot u)
   +
\nu(1) \cdot \bv \cdot \TTch(u) 
   +
\sum_{u}  \nu(\bv \cdot u_{(1)}) \cdot \bv \cdot \TTch(u_{(2)}) \\
  & = &
2 \bv \cdot \TTch(u) 
   +
(\av-\bv) \cdot \sum_{u}
      \kappa(u_{(1)}) \cdot \bv \cdot \TTch(u_{(2)}) \\
  & = &
2 \bv \cdot \TTch(u) 
   +
(\av-\bv) \cdot \sigma(u)   .
\end{eqnarray*}

For the operator 
$\sigma$ we have
\begin{eqnarray*}
\sigma(\av \cdot u)
  & = &
\kappa(1) \cdot \bv \cdot \TTch(u)
    +
\sum_{u} \kappa(\av \cdot u_{(1)}) \cdot \bv \cdot \TTch(u_{(2)})   \\
  & = &
\bv \cdot \TTch(u)
    +
(\av-\bv) \cdot \sum_{u} \kappa(u_{(1)}) \cdot \bv \cdot \TTch(u_{(2)})   \\
  & = &
\bv \cdot \TTch(u)
    +
(\av-\bv) \cdot \sigma(u) ,
\end{eqnarray*}
and
\begin{eqnarray*}
\hspace{30 mm}
\sigma(\bv \cdot u)
  & = &
\kappa(1) \cdot \bv \cdot \TTch(u)
    +
\sum_{u} \kappa(\bv \cdot u_{(1)}) \cdot \bv \cdot \TTch(u_{(2)})   \\
  & = &
\bv \cdot \TTch(u)    .
\hspace{55 mm}
\hspace{30 mm}
\qed
\end{eqnarray*}
\end{proof_qed}

\begin{corollary}
For an $\ab$-polynomial $u$ we have
$$\TTch((\av-\bv) \cdot u) = (\av-\bv) \cdot \TTch(u).$$
As a consequence, we have 
$$\TTch((\cv^{2}-2\dv) \cdot u) = (\cv^{2}-2\dv) \cdot \TTch(u).$$
\label{corollary_e}
\end{corollary}
\begin{proof}
The first part follows from subtracting equation~(\ref{equation_tau_b}) from
equation~(\ref{equation_tau_a}).
The second part follows from $(\av-\bv)^{2} = \cv^{2}-2\dv$.
\end{proof}

Define the operator $\pi$ on $\zab$ by
$$   \pi(u)
    =
2 \bv \cdot \TTch(u)
   +
2 (\av - \bv) \cdot \sigma(u) .  $$
We now restrict our attention to the subalgebra $\zcd$.

\begin{proposition}
The operator $\TTch$ satisfies the following joint recursion
with the operator $\pi$:
\begin{eqnarray}
\TTch(\av) & = & \cv \\
\TTch(\cv \cdot u)
  & = & 
\cv \cdot \TTch(u)
  +
\pi(u)  ,   \\
\TTch(\dv \cdot u)
  & = & 
2 \dv \cdot \TTch(u)
  +
\cv \cdot \pi(u)   , \label{equation_third}\\
\pi(\av)
  & = &
2 \dv  ,  \\
\pi(\cv \cdot u)
  & = & 
2\dv \cdot \TTch(u)
  +
\cv \cdot \pi(u)    ,   \\
\pi(\dv \cdot u)
  & = & 
\cv 2\dv \cdot \TTch(u)
  +
2\dv \cdot \pi(u)   . 
\end{eqnarray}
\label{proposition_cd_recursion}
\end{proposition}
\begin{proof}
By Proposition~\ref{proposition_ab_recursion}
we have
\begin{eqnarray}
\TTch(\cv \cdot u)
    & = & 
\TTch(\av \cdot u)
      +
\TTch(\bv \cdot u)   \nonumber  \\
    & = & 
(\av + \bv) \cdot \TTch(u)
   +
2 \bv \cdot \TTch(u) 
   +
2 (\av-\bv) \cdot \sigma(u)    \label{equation_two}\\
    & = & 
\cv \cdot \TTch(u)
  +
\pi(u) .  \nonumber
\end{eqnarray}
Iterating Proposition~\ref{proposition_ab_recursion}
twice yields
\begin{eqnarray}
\TTch(\dv \cdot u)
    & = & 
\TTch(\av \bv \cdot u)
      +
\TTch(\bv \av \cdot u)   \nonumber \\
    & = & 
(\av + \bv) \cdot \TTch(\bv \cdot u)
   +
(\av - \bv) \cdot \sigma(\bv \cdot u)
   +
2 \bv \cdot \TTch(\av \cdot u) 
   +
(\av-\bv) \cdot \sigma(\av \cdot u) \nonumber \\
    & = & 
\left[
2 \cv \bv 
  +
(\av - \bv) \cdot \bv
  +
2 \bv \cdot (\av + \bv) 
  +
(\av - \bv) \cdot \bv
\right] \cdot \TTch(u)
  \nonumber \\
& &
+
\left[
\cv \cdot (\av-\bv)
  +
2 \bv \cdot (\av - \bv) 
  +
(\av - \bv) \cdot (\av - \bv) 
\right] \cdot \sigma(u)  \nonumber \\
    & = & 
\left(
2 \dv
  +
2 \cv \cdot \bv 
\right) \cdot \TTch(u)
  +
2 \cv \cdot (\av-\bv) \cdot \sigma(u)  \label{equation_smile}\\
    & = & 
2 \dv \cdot \TTch(u)
  +
\cv \cdot \pi(u)    .  \nonumber
\end{eqnarray}

For the operator $\pi$ we have
\begin{eqnarray*}
\pi(\cv \cdot u)
  & = &
2 \bv \cdot \TTch(\cv \cdot u)
   +
2 (\av - \bv) \cdot \sigma(\cv \cdot u) \\
  & = &
2 \bv \cdot
\left[
(\av + \bv) \cdot \TTch(u)
   +
2 \bv \cdot \TTch(u) 
   +
2 (\av-\bv) \cdot \sigma(u)
\right]
\\
& &
   +  \;
2 (\av - \bv) \cdot
\left[
2 \bv \cdot \TTch(u)
    +
(\av-\bv) \cdot \sigma(u)
\right] \\
  & = &
(2\dv + 2 \cv\bv) \cdot \TTch(u)
  +
\cv \cdot 2 (\av-\bv) \cdot \sigma(u) \\
  & = &
2\dv \cdot \TTch(u)
  +
\cv \cdot \pi(u)   .
\end{eqnarray*}
Here for the second equality we have applied~(\ref{equation_two}).

A straightforward double iteration of 
Proposition~\ref{proposition_functional_equation}
yields
\begin{eqnarray}
\sigma(\dv \cdot u)
  & = &
\sigma(\av \bv \cdot u)
  +
\sigma(\bv \av \cdot u)    \nonumber \\
  & = &
\bv \cdot \TTch(\bv \cdot u)
    +
(\av-\bv) \cdot \sigma(\bv \cdot u) 
    +
\bv \cdot \TTch(\av \cdot u) \nonumber \\
  & = &
\left[
2 \bv^{2}
  +
(\av-\bv) \cdot \bv
  +
\bv \cdot (\av + \bv)
\right] \cdot \TTch(u)
    +
2 \bv \cdot (\av-\bv) \cdot \sigma(u) \nonumber \\
  & = &
(\dv
  +
2 \bv^{2}) \cdot \TTch(u)
    +
2 \bv \cdot (\av-\bv) \cdot \sigma(u)    . \label{equation_four}
\end{eqnarray}

Finally, we have
\begin{eqnarray*}
\pi(\dv \cdot u)
  & = &
2 \bv \cdot \TTch(\dv \cdot u)
   +
2 (\av - \bv) \cdot \sigma(\dv \cdot u) \\
  & = &
2 \bv \cdot
\left[
\left(
2 \dv
  +
2 \cv \cdot \bv 
\right) \cdot \TTch(u)
  +
2 \cv \cdot (\av-\bv) \cdot \sigma(u)
\right]
   \\
& & + \;
2 (\av - \bv) \cdot
\left[
(\dv
  +
2 \bv^{2}) \cdot \TTch(u)
    +
2 \bv \cdot (\av-\bv) \cdot \sigma(u)
\right] \\
  & = &
\left[
4 \bv \dv
  +
4 \bv \cv \bv
  +
2 (\av-\bv) \dv
  +
4 (\av - \bv) \bv^{2}
\right] \cdot \TTch(u)
\\
& &
  +  \;
\left[
4 \bv \cdot \cv \cdot (\av-\bv)
  +
4 (\av-\bv) \bv (\av-\bv)
\right] \cdot \sigma(u)   \\
  & = &
\left(
2 \cv \dv
  +
4 \dv \bv
\right) \cdot \TTch(u)
  +
2 \dv \cdot 2 (\av-\bv) \cdot \sigma(u)   \\
  & = &
\cv 2\dv \cdot \TTch(u)
  +
2\dv \cdot \pi(u)   ,
\end{eqnarray*}
where the second equality follows from~(\ref{equation_third}),
(\ref{equation_two}),
(\ref{equation_smile}) and~(\ref{equation_four})
\end{proof}

Note that
Proposition~\ref{proposition_cd_recursion}
offers different proofs for 
Theorem~\ref{theorem_non-negativity},
Corollaries~\ref{corollary_ctd}
and~\ref{corollary_2_power}.

The next proposition relates
the operators
$\TTch$ and $\pi$ with the operator $\omega$.
\begin{proposition}
For $\cd$-polynomials $v$ we have 
\begin{eqnarray*}
        \TTch(v \cdot \av)
  & = &
        \omega(\av \cdot v^{*})^{*}   ,          \\
        \pi(v \cdot \av)
  & = &
        \omega(\av \cdot v^{*} \cdot \bv)^{*}  .
\end{eqnarray*}
\end{proposition}
This proposition is straightforward to prove
using induction, and hence we omit the proof.
Notice this argument offers a second proof of
Theorem~\ref{theorem_TTch_omega}.

The next relation extends a result
from~\cite{Hetyei_Tchebyshev}
where the special case of the ladder poset was considered.
\begin{corollary}
For all $\ab$-polynomials $u$
we have 
$$   \TTch(\cv^{2} \cdot u)
  =
2 \cv \cdot \TTch(\cv \cdot u)
  +
(2\dv - \cv^{2}) \cdot \TTch(u)  .    $$
\label{corollary_cv_cv}
\end{corollary}
\begin{proof_qed}
{}From Proposition~\ref{proposition_cd_recursion}, we have
\begin{eqnarray*}
\hspace{30 mm}
\TTch(\cv^{2} \cdot u)
  & = & 
\cv \cdot \TTch(\cv \cdot u)
  +
\pi(\cv \cdot u) \\  
  & = & 
\cv \cdot \TTch(\cv \cdot u)
  +
2\dv \cdot \TTch(u)
  +
\cv \cdot \pi(u)    \\
  & = & 
\cv \cdot \TTch(\cv \cdot u)
  +
2\dv \cdot \TTch(u)
  +
\cv \cdot 
\left(
\TTch(\cv \cdot u)
   -
\cv \cdot \TTch(u)
\right)  \\
  & = & 
2 \cv \cdot \TTch(\cv \cdot u)
  +
(2\dv - \cv^{2}) \cdot \TTch(u)  .
\hspace{27 mm}
\hspace{30 mm}
\qed
\end{eqnarray*}
\end{proof_qed}

As a corollary to this recursion, we can now explain
the name Tchebyshev. This result is due to Hetyei, who
studied the Tchebyshev transform of the ladder poset.
Recall the ladder poset of rank $n+1$ is the unique poset
with $\cd$-index $\cv^{n}$.

\begin{corollary}
Substituting $\cv$ to be $x$ and $\dv$ to be $(x^{2}-1)/2$
in $\TTch(\cv^{n-1} \cdot \av)$ we obtain the Tchebyshev polynomial
of the first kind $T_{n}(x)$. 
\end{corollary}
Under this substitution the recurrence in
Corollary~\ref{corollary_cv_cv} becomes the
recurrence for the Tchebyshev polynomials.
It remains to observe that
the substitution takes
$\TTch(\av)$ and $\TTch(\cv \cdot \av)$
to
$T_{1}(x) = x$ and
$T_{2}(x) = 2x^{2} - 1$, respectively.

\section{$EL$-shellability}
\setcounter{equation}{0}

For a poset $P$ let ${\cal H}(P)$ be the set of
edges in the Hasse diagram of $P$, that is,
${\cal H}(P)
   =
      \{ (x,y) \:\: : \:\: x,y \in P, \: x \prec y \}$.
An {\em $R$-labeling} of a poset $P$ is a map $\lambda$ from
${\cal H}(P)$ to $\Lambda$,
a linearly ordered set of labels,  such that
in every interval $[x,y]$ there is a unique maximal
(saturated) chain
$x = x_{0} \prec x_{1} \prec \cdots \prec x_{k} = y$
having increasing labels,
that is,
$\lambda(x_{0},x_{1}) \leq_{\Lambda}
 \lambda(x_{1},x_{2}) \leq_{\Lambda} \cdots \leq_{\Lambda}
 \lambda(x_{k-1},x_{k})$.
Such a chain is called {\em rising}.
Furthermore an $R$-labeling is an {\em $EL$-labeling}
if the unique rising chain in every interval is
also the maximal chain with the lexicographically least labels.
A poset having an $EL$-labeling is said to be {\em $EL$-shellable}.
For further information regarding $EL$-labelings and their
topological consequences, see for example~\cite{Bjorner_Wachs}.

Recall the Jordan-H\"older set $\JH(x,y)$ of an interval $[x,y]$
is the collection of all strings of
labels occurring from the maximal chains in the interval,
that is,
$$   \JH(x,y)
   =
     \{ 
        (\lambda(x_{0},x_{1}),
         \lambda(x_{1},x_{2}), \ldots,
         \lambda(x_{k-1},x_{k}))
             \:\: : \:\:
         x = x_{0} \prec x_{1} \prec \cdots \prec x_{k} = y  \}  .  $$

\begin{theorem}
Let $P$ be an $EL$-shellable poset.
Then the Tchebyshev transform $T(P)$ is
an $EL$-shellable poset.
\label{theorem_EL_shellability}
\end{theorem}
\begin{proof}
Suppose the poset $P$ has label set
$\Lambda = \{\lambda_{1} < \cdots < \lambda_{k}\}$.
Define the new label set
$\Gamma = \{\lambda_{1}^s < \cdots < \lambda_{k}^s < 0 < 
            \lambda_{1}^b < \cdots < \lambda_{k}^b\}$.
Here one should think of the superscript $s$ as denoting
``small'' labels and the superscript $b$ as denoting
``big'' labels.
In 
the Tchebyshev poset $T(P)$ label
the edges in the Hasse diagram by the following rule:
$$
\left\{
\begin{array}{c c c}
\lambda([x,y],[y,w]) & = & \lambda(y,w)^{s} , \\
\lambda([x,y],[x,w]) & = & \lambda(y,w)^{b} , \\
\lambda([x,y],\ho_{T(P)}) & = & 0 .
\end{array}
\right.
$$
We claim this is an $EL$-labeling of the Tchebyshev poset
$T(P)$.
For a set $X$ of strings of labels from the set $\Lambda$,
let $X^{s}$ and $X^{b}$ denoted the set of strings where each
label has been signed with $s$, respectively $b$.
Similarly, let $X^{sb}$ denote
the set of strings where each
label has arbitrarily been signed $s$ or $b$.

There are three types of intervals to consider.
\begin{itemize}
\item[(i)]
An interval of the form
$I = [[x,y],[x,w]]$ in $T(P)$ is isomorphic to the interval
$[y,w]$ in the original poset~$P$.
In this case, 
the edge labels
are from the set
$\{ \lambda_1^b, \ldots, \lambda_k^b\}$
and the Jordan-H\"older set of the interval $I$
is described by
$$    \JH(y,w)^{b}   .   $$
Hence the lexicographically least maximal chain
in the interval $I$ is to take the
lexicographically least maximal chain
in the interval $[y,w]$ and change the
labels $\lambda_i$ to $\lambda_i^b$.

\item[(ii)] 
Let $I$ be an interval
of the form $[[x,y],[z,w]]$,
where 
$z$ is an element of rank $j$ from the half-open interval
$[y,w)$ in the poset $P$.
Observe that $0 \leq j < k$.
Any maximal chain
$\{ [x,y] = [x_{0},y_{0}] \prec [x_{1},y_{1}] \prec
       \cdots \prec [x_{k},y_{k}] = [z,w] \}$
in the interval $I$
satisfies
$\{ y = y_{0} \prec y_{1} \prec \cdots \prec y_{k} = w \}$
is a maximal chain in the interval $[y,w]$
and $z = y_{j} = x_{j+1} = \cdots = x_{k}$.
Thus the Jordan-H\"older set of the
interval $I$
is described by 
$$ \bigcup_{z \prec y_{j+1}}
      \JH(y,z)^{sb}  \circ
      \lambda(z,y_{j+1})^{s}  \circ
      \JH(y_{j+1},w)^{b}  , $$
where $\circ$ denotes concatenation.
To obtain a rising chain
in the interval $I$,
let
$m = \{ y = y_{0} \prec y_{1} \prec \cdots \prec y_{j} = z \}$
be the unique rising chain in
the interval $[x,z]$
and let
$m' = \{ z = y_{j} \prec y_{j+1} \prec \cdots \prec y_{k} = w \}$
be the unique rising chain in
the interval $[z,w]$.
Set $x_{i} = y_{i-1}$ for $0 < i \leq j$
and $x_{i} = z$ for $j+1 \leq i \leq k$.
The string of labels of this maximal chain is given by
$$     (\lambda(y_{0},y_{1})^{s},
       \lambda(y_{1},y_{2})^{s},
          \ldots,
       \lambda(y_{j-1},y_{j})^{s},
       \lambda(y_{j},y_{j+1})^{s},
       \lambda(y_{j+1},y_{j+2})^{b},
          \ldots,
       \lambda(y_{k-1},y_{k})^{b})        .   $$
It is straightforward to see that
this chain is the unique rising and lexicographic least
maximal chain in the interval.

\item[(iii)] 
Let $I$ be the interval
of the form $[[x,y],\ho_{T(P)}]$.
Any maximal chain
$\{ [x,y] = [x_{0},y_{0}] \prec [x_{1},y_{1}] \prec
       \cdots \prec [x_{k},y_{k}] \prec \ho_{T(P)} \}$
in the interval $I$
satisfies
$\{ y = y_{0} \prec y_{1} \prec \cdots \prec y_{k} = \ho_{P} \}$
is a maximal chain in the interval $[y,\ho_{P}]$.
Thus the Jordan-H\"older set of the
interval $I$
is described by 
$$    \JH(y,\ho_{P})^{sb}  \circ  0    . $$
Since all the labels signed with $s$ are smaller than $0$,
a rising chain can only have these ``small'' labels.
The unique
rising chain
in the interval $[y,\ho_{P}]$
is 
$\{ y = y_{0} \prec y_{1} \prec \cdots \prec y_{k} = \ho_{P} \}$.
To obtain the desired maximal chain in the interval $I$
with the correct labels,
let 
$x_{i} = y_{i-1}$
for $0 < i \leq k$.
This rising chain is also the lexicographic least.
\end{itemize}
Hence we conclude $T(P)$ has an $EL$-labeling.
\end{proof}

As a corollary to
Theorem~\ref{theorem_EL_shellability}
and its proof
we have:
\begin{corollary}
Let $P$ be a poset with an $R$-labeling having label set $\Lambda$.
Then the Tchebyshev transform $T(P)$ has
an $R$-labeling with 
the label set given by
$\Lambda^{sb} \cup \{0\}$
and
the Jordan-H\"older set
given by
$\JH(T(P)) = \JH(P)^{sb} \circ  0 $.
\label{corollary_R}
\end{corollary}

\section{The Tchebyshev transform of Cartesian products}
\label{section_Tchebyshev_Cartesian}
\setcounter{equation}{0}

In the papers~\cite{Ehrenborg_Fox,Ehrenborg_Readdy_c},
Ehrenborg-Fox and Ehrenborg-Readdy
studied 
the behavior of the $\cd$-index
under
the {\em Cartesian product} $P \times Q$ and
the {\em diamond product} $P \diamond Q$,
where $P$ and $Q$ are posets.
This latter product 
is defined as $P \diamond Q =
(P - \{\hz\}) \times (Q - \{\hz\}) \cup \{\hz\}$.
For our purposes, we need to consider the dual of the diamond product,
namely 
$$  P \diamond^{*} Q
     =
    (P - \{\ho\}) \times (Q - \{\ho\}) \cup \{\ho\}   .   $$
In other words,
$P \diamond^{*} Q = (P^{*} \diamond Q^{*})^{*}$.

We have the following result.
\begin{theorem}
Given two posets $P$ and $Q$,
the flag $f$-vector of the 
Tchebyshev transform of the Cartesian product $P \times Q$
is equal to the flag $f$-vector of the
dual diamond product of the two Tchebyshev transforms
$T(P)$ and $T(Q)$, that is,
$$  \Psi(T(P \times Q)) = \Psi(  T(P) \diamond^{*} T(Q)  )  .  $$
\label{theorem_Cartesian}
\end{theorem}

\begin{figure}
\setlength{\unitlength}{0.50mm}
\begin{center}
\begin{picture}(60,60)(0,0)
\put(10,0){\circle*{3}}
\put(30,10){\circle*{3}}
\put(50,0){\circle*{3}}
\put(10,60){\circle*{3}}
\put(30,50){\circle*{3}}
\put(50,60){\circle*{3}}

\put(10,0){\line(1,0){40}}
\put(10,0){\line(2,1){20}}
\put(30,10){\line(2,-1){20}}
\put(10,60){\line(1,0){40}}
\put(10,60){\line(2,-1){20}}
\put(30,50){\line(2,1){20}}

\qbezier(10,0)(5,30)(10,60)
\qbezier(10,0)(15,30)(10,60)
\qbezier(30,10)(25,30)(30,50)
\qbezier(30,10)(35,30)(30,50)
\qbezier(50,0)(45,30)(50,60)
\qbezier(50,0)(55,30)(50,60)

\end{picture}
\end{center}
\caption{A $CW$-decomposition of the $2$-sphere
with three $2$-gons, two triangles and three squares.
The face lattice is the Tchebyshev transform of $B_{3}$.}
\label{figure_picture}
\end{figure}
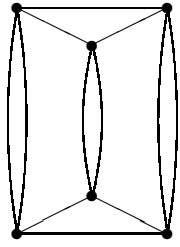

In general, it is not true that
the two posets
$T(P \times Q)$ and $T(P) \diamond^{*} T(Q)$
are isomorphic. A counterexample
is to take $P = B_{2}$ and $Q = B_{1}$.
The Tchebyshev transform of $B_{2}$ is isomorphic to
the face lattice of
a square and 
the Tchebyshev transform of $B_{1}$ is the face lattice of
a line segment.
Hence 
$T(B_{2}) \diamond^{*} T(B_{1})$
is the face lattice of the $3$-dimensional crosspolytope.
However the Tchebyshev transform of
$B_{2} \times B_{1} = B_{3}$ is not a lattice.
It is the face poset of the $CW$-complex displayed in
Figure~\ref{figure_picture}.

Observe that an alternate proof of
Corollary~\ref{corollary_Tchebyshev_the_Boolean_algebra}
follows directly from
Theorem~\ref{theorem_Cartesian} by considering
the Boolean algebra $B_{n} = B_{1}^{n}$.

To prove Theorem~\ref{theorem_Cartesian} we need
the following result
from~\cite{Ehrenborg_Readdy_c}:
\begin{theorem}
There exists two bilinear operators
$M$ and $N$ on $\zab$ such that
for two graded posets $P$ and $Q$ we have
\begin{eqnarray}
\Psi(P \times Q)
  & = &
M(\Psi(P), \Psi(Q)) , 
\label{equation_M} \\
\Psi(P \diamond Q)
  & = &
N(\Psi(P), \Psi(Q)) .
\label{equation_N} 
\end{eqnarray}
\label{theorem_M_and_N}
\end{theorem}
Recursions for the two bilinear operators
$M$ and $N$ have been developed in~\cite{Ehrenborg_Fox}.
Defining $N^{*}$ by
$N^{*}(u,v) = N(u^{*},v^{*})^{*}$, we have
\begin{equation}
\Psi(P \diamond^{*} Q)
    =
N^{*}(\Psi(P), \Psi(Q)) .
\label{equation_N_star} 
\end{equation}

Theorem~\ref{theorem_M_and_N}
states that on the flag $f$-vector level
the Cartesian product and the dual diamond product
are bilinear. Hence Theorem~\ref{theorem_Cartesian}
can be reformulated as follows.
\begin{theorem}
Given $\ab$-polynomials $u$ and $v$,
we have
$$     \TTch( M(u,v) \cdot \av )
   =
       N^{*}(
              \TTch( u \cdot \av ) ,
              \TTch( v \cdot \av ) 
             )                         .   $$
\label{theorem_ab_Cartesian}
\end{theorem}

Notice that from the results of Section~\ref{section_omega},
Theorem~\ref{theorem_Cartesian} is true for $\cd$-polynomials.

To prove Theorem~\ref{theorem_Cartesian}
it is enough to prove the identity for a class
of posets having $\ab$-indexes which span $\zab$.
We will prove the identity for 
posets that admit $R$-labelings.
\begin{proposition}
Let $P_{1}$ and $P_{2}$ be two posets
such that each has an $R$-labeling.
Then we have 
$$  \Psi(T(P_{1} \times P_{2})) 
  =
    \Psi(  T(P_{1}) \diamond^{*} T(P_{2})  )  .  $$
\label{proposition_R_Cartesian}
\end{proposition}

Let $P$ be a graded poset of rank $n+1$ that has an $R$-labeling. The
strings of labels in the Jordan-H\"older set $\JH(P)$ have length
$n+1$. For such a string $\lambda = (\lambda_{1}, \lambda_{2}, \ldots,
\lambda_{n+1})$, define its {\em decent word} to be
$u_{\lambda} = u_{1} u_{2} \cdots u_{n}$ by letting
$u_{i} = \av$ if $\lambda_{i} \leq \lambda_{i+1}$
and
$u_{i+1} = \bv$ otherwise.
Then we have the following result which expresses the $\ab$-index of
the poset $P$ in terms of the Jordan-H\"older set $\JH(P)$.
\begin{proposition}
Let $P$ be a poset with an $R$-labeling. Then the
$\ab$-index of $P$ is given by
$$    \Psi(P) = \sum_{\lambda \in \JH(P)} u_{\lambda}  .  $$
\label{proposition_Bjorner_Stanley}
\end{proposition}
The original formulation of this result is
due to Bj\"orner-Stanley~\cite{Bjorner}.
The reformulation 
in
Proposition~\ref{proposition_Bjorner_Stanley}
can be found
in~\cite{Billera_Ehrenborg_Readdy_om}.

Given two strings
$\xx = (x_{1}, \ldots, x_{n})$
and
$\yy = (y_{1}, \ldots, y_{m})$, define
their {\em shuffle product} $\xx \star \yy$
to be set of all ${{n+m} \choose n}$ shuffles of them,
that is,
\begin{eqnarray*}
  \xx \star \yy
  & = &
    \{ (z_{1}, \ldots, z_{n+m}) \:\: : \:\:
       z_{i_{p}} = x_{p}, 
       z_{j_{q}} = y_{q}, \\
  &   &
\mbox{ where }
          \{i_{1} < \cdots < i_{n}\}
        \cup
          \{j_{1} < \cdots < j_{m}\}
        =
          \{1, \ldots, n+m\}     \}  .  
\end{eqnarray*}
For two sets of strings $X$ and $Y$, define their shuffle
product be
$$   X \star Y
        =
     \bigcup_{\xx \in X, \: \yy \in Y}   \xx \star \yy  .  $$

\begin{lemma}
For $i = 1,2$ let $P_{i}$ be a poset of rank $n_{i}$
with an $R$-labeling 
$\lambda^{i}$ and linearly ordered label poset $\Lambda^{i}$.
Without loss of generality assume that 
$\Lambda^{1}$ and $\Lambda^{2}$ are disjoint.
Let $\Gamma$ be a linear extension of the union
$\Lambda^{1} \cup \Lambda^{2}$.
Then $P_{1} \times P_{2}$ has an $R$-labeling $\gamma$
given by
$$      \left\{ \begin{array}{c c c}
           \gamma((x,y),(z,y)) & = & \lambda^{1}(x,z) , \\
           \gamma((x,y),(x,w)) & = & \lambda^{2}(y,w) . \\
                \end{array} \right.   $$
Moreover, the Jordan-H\"older set
$\JH(P_{1} \times P_{2})$ is given by all shuffle products
of the strings from
$\JH(P_{1})$ and $\JH(P_{2})$,
that is,
$$ \JH(P_{1} \times P_{2})
       =
   \JH(P_{1}) \star \JH(P_{2})   .  $$
\label{lemma_JH_Cartesian}
\end{lemma}
\begin{proof}
Every maximal chain in the product $P_{1} \times P_{2}$
comes from one maximal chain in $P_{1}$ and one maximal chain in
$P_{2}$.
Conversely, for each pair $(c_{1}, c_{2})$ of maximal chains, where $c_{i}$
is a maximal chain in $P_{i}$,
there are ${{n_{1} + n_{2}} \choose n_{1}}$
maximal chains in $P_{1} \times P_{2}$.
Moreover, the labels of these ${{n_{1} + n_{2}} \choose n_{1}}$
maximal chains are the shuffle product of the labels of $c_{1}$
and the labels of $c_{2}$. Hence the Jordan-H\"older set
of the Cartesian product $P_{1} \times P_{2}$ has the desired
form.

Consider an interval $I = [(x_{1},x_{2}), (y_{1},y_{2})]$
in the product $P_{1} \times P_{2}$.
Let $m^{i}$
be the string of labels of a maximal chain
in the interval $[x_{i}, y_{i}]$ in the poset $P_{i}$.
If $m^{1}$ or $m^{2}$ 
has a descent then
all the strings of labels in the shuffle product
$m^{1} \star m^{2}$ have at least one descent.
Now let $\lambda^{i}$ be the string of labels of the unique rising chain
in the interval $[x_{i}, y_{i}]$.
Then there is exactly one shuffle among $m^{1} \star m^{2}$
that is a rising string. 
Hence the interval $I$
has a unique rising chain, proving 
$\mu$ is an $R$-labeling.
\end{proof}

\begin{lemma}
For $i = 1,2$ let $P_{i}$ be a poset with an $R$-labeling 
$\lambda^{i}$ and linearly ordered label poset $\Lambda^{i}$.
Assume that each edge in the Hasse diagram 
between a coatom of $P_{i}$
and
the maximal element $\ho_{P_{i}}$ is labeled $0$ and no other
labels are equal to $0$. This condition can be expressed
as
$$ \lambda^{i}(x,y) = 0
       \:\: \Longleftrightarrow \:\:
   x \prec y = \ho_{P_{i}}  .  $$
Without loss of generality assume that 
$\Lambda^{1} \cap \Lambda^{2} = \{0\}$.
Let $\Gamma$ be a linear extension of the union
$\Lambda^{1} \cup \Lambda^{2}$.
Then $P_{1} \diamond^{*} P_{2}$ has an $R$-labeling $\gamma$
given by
$$      \left\{ \begin{array}{c c c}
           \gamma((x,y),(z,y)) & = & \lambda^{1}(x,z) , \\
           \gamma((x,y),(x,w)) & = & \lambda^{2}(y,w) , \\
           \gamma((x,y),\ho_{P_{1} \diamond^{*} P_{2}})
                              & = & 0 . \\
                \end{array} \right.   $$
Moreover, the Jordan-H\"older set is given by
$$  \JH(P_{1} \diamond^{*} P_{2})
   =
    \left(
      \JH_{0}(P_{1}) \star \JH_{0}(P_{2})
    \right) \circ 0   ,  $$
where $\JH_{0}(P_{i})$ is the set of all the strings
in the Jordan H\"older set 
$\JH(P_{i})$ with the $0$ at the end removed.
\label{lemma_JH_diamond}
\end{lemma}
\begin{proof}
Directly from the identity
$(P_{1} - \{\ho_{P_{1}}\}) \times (P_{2} - \{\ho_{P_{2}}\})
  =
(P_{1} \diamond^{*} P_{2}) - \{ \ho_{P_{1} \diamond^{*} P_{2}} \}$
it follows that
$\JH_{0}(P_{1} \diamond^{*} P_{2})
   =
 \JH_{0}(P_{1}) \star \JH_{0}(P_{2})$,
thus verifying the Jordan H\"older set
of the dual diamond product is as
described.

It remains to observe that $\gamma$ is an $R$-labeling.
By the same reasoning as in the proof of
Lemma~\ref{lemma_JH_Cartesian},
each interval of the form
$[(x,y),(z,w)]$ has a unique rising chain.
Hence it is enough
to show each interval of the form
$I = [(x,y), \ho_{P_{1} \diamond^{*} P_{2}}]$
has a unique rising chain.

Let $m^{i} \circ 0$
be the string of labels of a maximal chain
in the interval $[x_{i}, \ho_{P_{i}}]$ in the poset $P_{i}$.
If $m^{1} \circ 0$ or $m^{2} \circ 0$ 
has a descent then
all the strings of labels in the shuffle product
$(m^{1} \star m^{2}) \circ 0$ has at least one descent.
Now let $m^{i} \circ 0$
be the string of labels of the unique rising chain
in the interval $[x_{i}, \ho_{P_{i}}]$.
Then there is exactly one shuffle among
$(\lambda^{1} \star \lambda^{2}) \circ \{0\}$
that is a rising string. 
Hence the interval $I$
has a unique rising chain, proving 
$\gamma$ is an $R$-labeling.
\end{proof}

\begin{proof_}{Proof of Proposition~\ref{proposition_R_Cartesian}:}
Let the $R$-labeling of the poset
$P_{i}$ have label set $\Lambda_{i}$,
where we assume $\Lambda_{1}$ and $\Lambda_{2}$ are disjoint.
Then the Cartesian product $P_{1} \times P_{2}$ has an $R$-labeling
with the label set
$\Lambda_{1} \cup \Lambda_{2}$
and the Jordan-H\"older set
$\JH(P_{1} \times P_{2}) = \JH(P_{1}) \star \JH(P_{2})$.
Now by Corollary~\ref{corollary_R}
the Tchebyshev transform of the product $P_{1} \times P_{2}$
has an $R$-labeling with
label set
$(\Lambda_{1} \cup \Lambda_{2})^{sb} \cup \{0\}$
and Jordan-H\"older set
$(\JH(P_{1}) \star \JH(P_{2}))^{sb} \circ \{0\}$.

Similarly, by Corollary~\ref{corollary_R}
the Tchebyshev transform of the poset $P_{i}$
has
an $R$-labeling with label set $\Lambda_{i}^{sb} \cup \{0\}$
and Jordan-H\"older set
$\JH(P_{i})^{sb} \circ \{0\}$.
Now by Lemma~\ref{lemma_JH_diamond}
the diamond product $T(P_{1}) \diamond^{*} T(P_{2})$
has an $R$-labeling with
label set
$\Lambda_{1}^{sb} \cup \Lambda_{2}^{sb} \cup \{0\}$
and Jordan H\"older set
$(\JH(P_{1})^{sb} \star \JH(P_{2})^{sb}) \circ \{0\}$.

As sets, the two label sets agree:
$$ (\Lambda_{1} \cup \Lambda_{2})^{sb} \cup \{0\}
  =
   \Lambda_{1}^{sb} \cup \Lambda_{2}^{sb} \cup \{0\}    .   $$
Additionally, as linearly ordered sets they are also
equal, since we can first choose the linear extension
of $\Lambda_{1} \cup \Lambda_{2}$ to be the unique
linear order where all the labels from $\Lambda_{1}$
is an initial segment.
Moreover, choose
the linear extension of
$(\Lambda_{1}^{sb} \cup \{0\}) \cup (\Lambda_{2}^{sb} \cup \{0\})$
to be as 
$\Lambda_{1}^{s}, \Lambda_{2}^{s}, \{0\}, \Lambda_{1}^{b}, \Lambda_{2}^{b}$,
in that order.

Finally, observe that the Jordan-H\"older sets of
the two posets
$T(P_{1} \times P_{2})$ and $T(P_{1}) \diamond^{*} T(P_{2})$
also are equal, namely,
$$   (\JH(P_{1}) \star \JH(P_{2}))^{sb} \circ \{0\}
   =
     \JH(P_{i})^{sb} \circ \{0\}         .  $$
Hence by Proposition~\ref{proposition_Bjorner_Stanley}
the posets have the same $\ab$-index.
\end{proof_}

\section{The Tchebyshev operator of the second kind}
\label{section_Tchebyshev_operator}
\setcounter{equation}{0}

Following Hetyei,
we will now define the Tchebyshev operator of the
second kind.
We demonstrate it is an algebra map with respect to
the mixing operator $M$ and a coalgebra map
with respect to the coproduct $\Delta$.
Moreover, we 
find the spectrum of this operator,
generalizing work in~\cite{Billera_Hsiao_van_Willigenburg}.

Define the two linear maps
$H, H^{*} : \zab \longrightarrow \zab$
by
$H(1) = H^{*}(1) = 0$
and
$H(\av \cdot u) = H(\bv \cdot u) = 
 H^{*}(u \cdot \av) = H^{*}(u \cdot \bv) = u$.
The map $H$ appears
in~\cite{Billera_Ehrenborg_Readdy_om}.
We have the following result from the same reference.
\begin{lemma}
For a poset $P$ of rank at least $2$ we have
\begin{eqnarray*}
  H(\Psi(P))
     & = &
   \sum_{a} \Psi([a,\ho]), \\
  H^{*}(\Psi(P))
     & = &
   \sum_{c} \Psi([\hz,c]),
\end{eqnarray*}
where the first sum 
ranges over all atoms $a$ of the poset $P$
and the second sum
ranges over all coatoms $c$ of the poset $P$.
\end{lemma}
Observe both $H$ and $H^{*}$ restrict to $\zcd$ by
$H(\cv \cdot u) = 
 H^{*}(u \cdot \cv) = 2 u$,
$H(\dv \cdot u) = \cv \cdot u$
and
$H^{*}(u \cdot \dv) = u \cdot \cv$.

\begin{definition}
The {\em Tchebyshev transform of the second kind}
is the linear map $\UTch : \zab \longrightarrow \zab$
defined by
$$    \UTch(u) =  H^{*}( \TTch( u \cdot \av) )   .  $$
\end{definition}
The explanation for this name is given by the next corollary.
This result is originally due to Hetyei.

\begin{corollary}
Substituting $\cv$ to be $x$ and $\dv$ to be $(x^{2}-1)/2$
in $1/2 \cdot \UTch(\cv^{n})$ yields the Tchebyshev polynomial
of the second kind $U_{n}(x)$. 
\end{corollary}
\begin{proof}
First observe that
under this substitution the
expressions
$1/2 \cdot \UTch(1) = 1$
and
$1/2 \cdot \UTch(\cv) = 2 \cv$
become
$U_{0}(x) = 1$ and $U_{1}(x) = 2x$.
Second, the recursion in 
Corollary~\ref{corollary_cv_cv}
transforms into
$\UTch(\cv^{2} \cdot u)
  =
 2 \cv \cdot \UTch(\cv \cdot u)
  +
(2\dv - \cv^{2}) \cdot \UTch(u)$.
Under the given substitution
this becomes the 
recursion for the
Tchebyshev polynomials of the second kind.
\end{proof}

\begin{proposition}
The Tchebyshev transform of the second kind
has the following expression:
$$
  \UTch(u) = \sum_{k \geq 1}
             \sum_{u}
             \nu(u_{(1)}) \cdot \bv \cdot
             \nu(u_{(2)}) \cdot \bv \cdots
                               \bv \cdot
             \nu(u_{(k)})   .
$$
\label{proposition_UTch}
\end{proposition}
\begin{proof}
By applying the definition 
of the Tchebyshev transform
appearing in
equation~(\ref{equation_TTch}),
we have
\begin{eqnarray*}
  \TTch(u \cdot \av)
  & = &
 \sum_{k \geq 1}
             \sum_{u}
             \nu(u_{(1)}) \cdot \bv \cdot
             \nu(u_{(2)}) \cdot \bv \cdots
                               \bv \cdot
             \nu(u_{(k-1)}) \cdot \bv \cdot
             \kappa(u_{(k)} \cdot \av)   \\
  &   &
    +
 \sum_{k \geq 2}
             \sum_{u}
             \nu(u_{(1)}) \cdot \bv \cdot
             \nu(u_{(2)}) \cdot \bv \cdots
                               \bv \cdot
             \nu(u_{(k-1)}) \cdot \bv \cdot
             \kappa(1)   \\
  & = &
             T(u) \cdot (\av-\bv) 
    +
 \sum_{k \geq 1}
             \sum_{u}
             \nu(u_{(1)}) \cdot \bv \cdot
             \nu(u_{(2)}) \cdot \bv \cdots
                               \bv \cdot
             \nu(u_{(k)}) \cdot \bv      .
\end{eqnarray*}
The result now follows by applying the map $H^{*}$.
\end{proof}

\begin{corollary}
The Tchebyshev transform of the second kind
is invariant under duality, that is,
$\UTch(u^{*}) = \UTch(u)^{*}$.
\end{corollary}

\begin{theorem}
The Tchebyshev transform of the second kind
is a coalgebra homomorphism, that is,
$$
  \Delta(\UTch(u)) = \sum_{u} \UTch(u_{(1)}) \tensor \UTch(u_{(2)})  .
$$
\label{theorem_UTch_coalgebra}
\end{theorem}
\begin{proof_qed}
Recall that
$\Delta(\nu(u)) = 0$.
By applying Proposition~\ref{proposition_UTch}, we obtain
\begin{eqnarray*}
\hspace*{15 mm}
\Delta(\UTch(u))
  & = &
             \sum_{k \geq 1}
             \sum_{u}
               \sum_{i=1}^{k-1}
             \nu(u_{(1)}) \cdot \bv \cdots \bv \cdot \nu(u_{(i)})
           \tensor
             \nu(u_{(i+1)}) \cdot \bv \cdots \bv \cdot \nu(u_{(k)})  \\
  & = &
   \sum_{i,j \geq 1}
      \sum_{u}
      \sum_{u_{(1)}}
      \sum_{u_{(2)}}
             \nu(u_{(1,1)}) \cdot \bv \cdots \cdot \nu(u_{(1,i)}) 
          \tensor
             \nu(u_{(2,1)}) \cdot \bv \cdots \cdot \nu(u_{(2,j)})   \\
  & = &
      \sum_{u}
           \UTch(u_{(1)}) \tensor \UTch(u_{(2)})  .
\hspace*{75 mm}
\hspace*{15 mm}
\qed
\end{eqnarray*}
\end{proof_qed}

\begin{theorem}
For two $\ab$-polynomials $u$ and $v$ we have
\begin{equation}
\UTch(M(u,v)) = M(\UTch(u),\UTch(v))  .
\end{equation}
In other words, 
the Tchebyshev transform of the second kind
is an algebra map under the product $M$.
\label{theorem_UTch_algebra}
\end{theorem}
\begin{proof_qed}
By Lemma~2.3 in~\cite{Ehrenborg_r-Birkhoff}
we have 
$$     H^{*}(N^{*}(u, v))
     =
       M(H^{*}(u), H^{*}(v))   .  $$
Applying $H^{*}$ to Theorem~\ref{theorem_ab_Cartesian},
we obtain
\begin{eqnarray*}
\hspace*{40 mm}
\UTch(u)
  & = &
H^{*}(\TTch( M(u,v) \cdot \av ) ) \\
  & = &
H^{*}(
       N^{*}(
              \TTch( u \cdot \av ) ,
              \TTch( v \cdot \av ) 
             )    ) \\
  & = &
M(
       H^{*}(
              \TTch( u \cdot \av ) ),
       H^{*}(
              \TTch( v \cdot \av ) )
  )  \\
  & = &
M(
              \UTch( u ),
              \UTch( v )
  )  .
\hspace*{20 mm}
\hspace*{40 mm}
\qed
\end{eqnarray*}
\end{proof_qed}

\begin{proposition}
Assume $u_{i}$ is an eigenvector with eigenvalue $\lambda_{i}$
of the Tchebyshev transform of the second kind $\UTch$
for $i=1,2$. Then
$M(u_{1},u_{2})$ is an eigenvector with
eigenvalue $\lambda_{1} \cdot \lambda_{2}$.
\label{proposition_M_eigenvalues}
\end{proposition}
\begin{proof}
Directly
$\UTch(M(u_{1},u_{2}))
    =
M(\UTch(u_{1}),\UTch(u_{2}))
    =
M(\lambda_{1} \cdot u_{1}, \lambda_{2} \cdot u_{2})
    =
\lambda_{1} \cdot \lambda_{2} \cdot M(u_{1}, u_{2})$.
\end{proof}

\begin{proposition}
Assume  $u$ is an eigenvector with eigenvalue $\lambda$
of the Tchebyshev transform of the second kind $\UTch$.
Then
$(\av-\bv) \cdot u$ is an eigenvector with
eigenvalue $\lambda$.
\label{proposition_a-b_eigenvalues}
\end{proposition}
\begin{proof}
Observe that
$\UTch((\av-\bv) \cdot u)
    =
H^{*}(\TTch((\av-\bv) \cdot u \cdot \av))
    =
(\av-\bv) \cdot H^{*}(\TTch(u \cdot \av))
    =
(\av-\bv) \cdot \UTch(u)$,
where the second step is by Corollary~\ref{corollary_e}.
\end{proof}

Let $\kab_{n}$ denote the set of all homogeneous $\ab$-polynomials
of degree $n$ with coefficients in the field $\kk$.
Hence the dimension of $\kab_{n}$
is $2^{n}$ and $\UTch_{n}$ is an endomorphism on $\kab_{n}$.

\begin{theorem}
Let $\UTch_{n}$ denote the restriction of $\UTch$ to $\ab$-polynomials
of degree $n$, that is, $\kab_{n}$.
Then the linear operator $\UTch_{n}$ is diagonalizable
and has the eigenvalue $2^{i+1}$ of multiplicity ${n \choose i}$
for $0 \leq i \leq n$.
Furthermore, a complete set of eigenvectors can be obtained
by starting with $1$ and repeatedly applying the two operations
\begin{eqnarray*}
      u & \longmapsto & \Pyr(u) = M(u,1)  , \\
      u & \longmapsto & L(u) = (\av-\bv) \cdot u , 
\end{eqnarray*}
$n$ times.
\label{theorem_UTch_spectrum}
\end{theorem}
\begin{proof}
Observe that
$1$ is an eigenvector with eigenvalue $2$.
By iterating Propositions~\ref{proposition_M_eigenvalues}
and~\ref{proposition_a-b_eigenvalues} $n$ times,
we obtain $2^{n}$ eigenvectors of degree $n$.
By Proposition~3.4 in~\cite{Billera_Hsiao_van_Willigenburg}
we know 
$$    \kab_{n+1} = \Pyr(\kab_{n}) \oplus L(\kab_{n})  .  $$
Hence
this set of eigenvectors is a complete set of eigenvectors,
that is, there are no linear dependencies among them.

Also since the pyramid operation $\Pyr$
multiplies an eigenvalue by $2$ and
the second operation $L$ preserves the eigenvalue, 
we may conclude the distribution of the eigenvalues of $\UTch_{n}$ is
precisely the binomial distribution.
\end{proof}

\section{A Hopf-algebra endomorphism on quasisymmetric functions}
\label{section_endomorphism}
\setcounter{equation}{0}

The main result of this section is prove 
the Tchebyshev transform of the second kind
is a Hopf algebra endomorphism.

Define the map $\UTch$ on a quasisymmetric function $f$
(where we intentionally use the same symbol
as the Tchebyshev transform of the second kind)
by
$$    \UTch(f) = \gamma(\UTch(\gamma^{-1}(f))), $$
where $f \in \QSym$ does not have a constant term.
Extend linearly to all quasisymmetric functions
by setting $\UTch(1_{\QSym}) = 1_{\QSym}$.

Theorems~\ref{theorem_UTch_algebra}
and~\ref{theorem_UTch_coalgebra}
imply the following result.
\begin{theorem}
The map $\UTch$ is a Hopf algebra
endomorphism on the Hopf algebra
of quasisymmetric functions.
\label{theorem_UTch_endomorphism}
\end{theorem}
\begin{proof_}{{\bf Sketch of proof:}}
We leave it to the reader to verify that $\UTch$
behaves well with the unit and the counit
of quasisymmetric functions.
Since the mixing operator $M$ on $\zab$
corresponds to the Cartesian product
on graded posets 
(equation~(\ref{equation_M}))
and the 
Cartesian product corresponds to
the product of quasisymmetric functions
(equation~(\ref{equation_product})),
it follows that
Theorem~\ref{theorem_UTch_algebra}
implies $\UTch$ is
algebra endomorphism on the quasisymmetric
functions.

Now for a quasisymmetric polynomial $f = \gamma(v)$, we have
\begin{eqnarray*}
\Delta^{\QSym}(\UTch(f))
  & = &
\Delta^{\QSym}(\gamma(\UTch(v)))  \\
  & = &
\gamma(\UTch(v)) \tensor 1_{\QSym}
  +
1_{\QSym} \tensor \gamma(\UTch(v))
  +
\sum_{v}
  \gamma(\UTch(v_{(1)})) \tensor \gamma(\UTch(v_{(2)})) \\
  & = &
\UTch(f) \tensor 1_{\QSym}
  +
1_{\QSym} \tensor \UTch(f)
  +
\sum_{v}
  \UTch(\gamma(v_{(1)})) \tensor \UTch(\gamma(v_{(2)})) \\
  & = &
(\UTch \tensor \UTch) \circ \Delta^{\QSym}(f)   ,
\end{eqnarray*}
where the second step is
that $\UTch$ is a coalgebra endomorphism
on $\zab$. This completes
the proof that~$\UTch$ is a coalgebra endomorphism on
quasisymmetric functions.
\end{proof_}

\section{Chain maps of the first and second kind}
\label{section_chain_maps}
\setcounter{equation}{0}

The results in
Sections~\ref{section_Tchebyshev_transform},
\ref{section_Tchebyshev_operator}
and~\ref{section_endomorphism}
motivate us to 
consider two general classes of maps.
In this section, we show one such class, the chain map of the second kind,
is a Hopf algebra endomorphism of quasisymmetric functions.

\begin{definition}
A character $G$ on $\kab$ is a functional
$G : \kab \longrightarrow \kk$
which is multiplicative with
respect to Cartesian product of posets,
that is,
$$   G(\Psi(P \times Q))
   =
     G(\Psi(P))
   \cdot
     G(\Psi(Q))   ,  $$
for all posets $P$ and $Q$ of rank greater than or equal to $1$.
\end{definition}

Theorem~\ref{theorem_UTch_endomorphism}
can be extended in the following manner.
Let $G$ be a character on $\kab$. Define
the functions $\widehat{g}$, $\widetilde{g}$ and $g$ on $\kab$ by
\begin{eqnarray*}
      \widehat{g}(u)
  & = &
      G(u) \cdot (\av-\bv)^{\deg(u)}  ,  \\
      \widetilde{g}(u)
  & = &
       \sum_{k \geq 1}
       \sum_{u}
           \widehat{g}(u_{(1)})
         \cdot
           \bv
         \cdot
           \widehat{g}(u_{(2)})
         \cdot
           \bv
         \cdots
           \bv
         \cdot
           \widehat{g}(u_{(k)})   , \\
      g(u)
  & = &
       \sum_{k \geq 1}
       \sum_{u}
           \kappa(u_{(1)})
         \cdot
           \bv
         \cdot
           \widehat{g}(u_{(2)})
         \cdot
           \bv
         \cdots
           \bv
         \cdot
           \widehat{g}(u_{(k)})   .  
\end{eqnarray*}
We call the maps $g$ and $\tilde{g}$, respectively,
{\em the chain maps of the first and second kind}.

\begin{examples}
{\rm
(i) $G$ always takes the value $1$.
Then $\widehat{g} = \kappa$ and the two maps $g$ and $\widetilde{g}$
are both equal to the identity map. 

\noindent
(ii) $G(\Psi(P)) = \sum_{x \in P} (-1)^{\rho(x)} \cdot \mu(\hz,x)$.
Then $g$ is the $\varphi$ map of oriented matroids
(see equation~(\ref{equation_phi}))
and
$\widetilde{g}$ is the Stembridge $\vartheta$ map.

\noindent
(iii)
An extension of the previous example is to take
$G(\Psi(P)) = \sum_{x \in P} (1-r)^{\rho(x)} \cdot \mu(\hz,x)$.
In this case, 
$g$ corresponds to $\varphi_{r}$ of the $r$-signed Birkhoff transform
and
$\widetilde{g}$ is the $r$-signed analogue of the
Stembridge map, $\vartheta_{r}$.

\noindent
(iv)
$G(\Psi(P))$ is the cardinality of the poset $P$.
In this case we have 
$g(u^{*})^{*}$ is the Tchebyshev transform of the first
kind
and
$\widetilde{g}(u^{*})^{*}$ is the Tchebyshev transform of the 
second kind.
}
\label{examples_four}
\end{examples}

\begin{proposition}
The following relations hold between the functions
$\widetilde{g}$ and $g$:
\begin{eqnarray}
\Delta(\widetilde{g}(u))
  & = &
\sum_{u} \widetilde{g}(u_{(1)}) \tensor \widetilde{g}(u_{(2)})  ,  
\label{equation_Delta_widetilde_g} \\
\Delta(g(u))
  & = &
\sum_{u} g(u_{(1)}) \tensor \widetilde{g}(u_{(2)})  , 
\label{equation_Delta_g} \\
g(\av \cdot u)
  & = &
       (\av-\bv) \cdot g(u)   
     +
       \bv \cdot \widetilde{g}(u)   , 
\label{equation_g_a_u}\\
g(\bv \cdot u)
  & = &
       \bv \cdot \widetilde{g}(u)   
\label{equation_g_b_u} .
\end{eqnarray}
\label{proposition_g_widetilde_g}
\end{proposition}
\begin{proof}
The proof that $\widetilde{g}$ is a coalgebra endomorphism
follows exactly along the same lines as the
proofs of Theorems~\ref{theorem_UTch_coalgebra}
and~\ref{theorem_UTch_endomorphism}.
The same proof idea also establishes
equation~(\ref{equation_Delta_g}).
Identity~(\ref{equation_g_a_u}) follows from
\begin{eqnarray*}
g(\av \cdot u)
  & = &
       \sum_{k \geq 1}
       \sum_{u}
           \kappa(\av \cdot u_{(1)})
         \cdot
           \bv
         \cdot
           \widehat{g}(u_{(2)})
         \cdot
           \bv
         \cdots
           \bv
         \cdot
           \widehat{g}(u_{(k)})    \\
   &   &
     +
       \sum_{k \geq 1}
       \sum_{u}
           \kappa(1)
         \cdot
           \bv
         \cdot
           \widehat{g}(u_{(1)})
         \cdot
           \bv
         \cdots
           \bv
         \cdot
           \widehat{g}(u_{(k)})    \\
  & = &
       (\av-\bv) \cdot g(u)   
     +
       \bv \cdot \widetilde{g}(u)   .
\end{eqnarray*}
Identity~(\ref{equation_g_b_u}) follows in a similar manner.
\end{proof}

\begin{proposition}
The chain map of the second kind $\widetilde{g}$
has the following form when applied to
the $\ab$-index, respectively, the quasisymmetric function
of a poset $P$:
\begin{eqnarray}
      \widetilde{g}(\Psi(P))
  & = &
       \sum_{c}
           G(\Psi([x_{0},x_{1}]))
         \cdots
           G(\Psi([x_{k-1},x_{k}]))
         \cdot
           \wt(c) ,
\label{equation_widetilde_g_Psi}                       \\[2 mm]
      \widetilde{g}(F(P))
  & = &
       \sum_{c}
           G(\Psi([x_{0},x_{1}]))
         \cdots
           G(\Psi([x_{k-1},x_{k}]))
         \cdot
           M_{(\rho(x_{0},x_{1}), \ldots, \rho(x_{k-1},x_{k}))}  ,
\label{equation_widetilde_g_F_P_chain}    \\[2 mm]
      \widetilde{g}(F(P))
  & = &
       \lim_{j \longrightarrow \infty}
       \sum_{m}
           G(\Psi([x_{0},x_{1}]))
         \cdots
           G(\Psi([x_{j-1},x_{j}]))
         \cdot
           t_{1}^{\rho(x_{0},x_{1})}
         \cdots
           t_{j}^{\rho(x_{j-1},x_{j})}    ,  
\label{equation_widetilde_g_F_P}
\end{eqnarray}
where the first two sums
are over all chains
$c = \{\hz = x_0 < x_1 < \cdots < x_k = \ho\}$
in the poset $P$ and the third sum is over all
multichains 
$m = \{\hz = x_{0} \leq x_{1} \leq \cdots \leq x_{j} = \ho\}$
in $P$.
\label{proposition_widetilde_g_Psi_F}
\end{proposition}
\begin{proof}
The first identity follows
by using the definition of $\widetilde{g}$
and the fact that the $\ab$-index is a coalgebra
homomorphism.
The second identity follows
from the first by applying the map $\gamma$.

To prove the third identity,
let $t_{j+1} = t_{j+2} = \cdots = 0$ in 
the second identity (identity~(\ref{equation_widetilde_g_F_P_chain})).
This 
restricts the sum to chains having at most $j$ steps,
that is, $k \leq j$.
Such chains can be expressed in terms of multichains with $j$ steps.
We do this by extending the composition $G \circ \Psi$ by
letting $G(\Psi(\bullet)) = 1$.  By the definition
of the monomial quasisymmetric function, we then have
\begin{eqnarray*}
  &   &
      \widetilde{g}(F(P)) |_{t_{j+1} = t_{j+2} = \cdots = 0} \\
  & = &
       \sum_{\hz = x_{0} \leq x_{1} \leq \cdots \leq x_{j} = \ho}
           G(\Psi([x_{0},x_{1}]))
         \cdots
           G(\Psi([x_{j-1},x_{j}]))   
         \cdot
           t_{1}^{\rho(x_{0},x_{1})}
         \cdots
           t_{j}^{\rho(x_{j-1},x_{j})}    .  
\end{eqnarray*}
Letting $j$ tend to infinity yields the desired identity.
\end{proof}

Observe in the proofs of
Propositions~\ref{proposition_g_widetilde_g}
and~\ref{proposition_widetilde_g_Psi_F}
we only used the fact that $G$ is functional on $\kab$,
not that $G$ is a character on $\kab$.

\begin{theorem}
Let $G$ be a character on $\kab$.
Then for all $\ab$-polynomials $u$ and $v$ we have
\begin{eqnarray*}
      \widetilde{g}(M(u,v))
  & = &
      M(\widetilde{g}(u),
        \widetilde{g}(v))    . 
\end{eqnarray*}
Equivalently, for all quasisymmetric functions $f_{1}$ and $f_{2}$
we have
\begin{eqnarray*}
      \widetilde{g}(f_{1} \cdot f_{2})
  & = &
      \widetilde{g}(f_{1}) \cdot \widetilde{g}(f_{2})  .
\end{eqnarray*}
\label{theorem_widetilde_M}
\end{theorem}
\begin{proof}
A multichain of length $m$ in the Cartesian product $P \times Q$
corresponds to two multichains
of length $m$,
with one coming
from  the poset $P$
and
the other from the poset $Q$.
By applying
equation~(\ref{equation_widetilde_g_F_P})
three times, we have
\begin{eqnarray*}
  &   &
      \widetilde{g}(F(P \times Q))         \\
  & = &
       \lim_{m \longrightarrow \infty}
       \sum_{\hz = (x_{0},y_{0}) \leq (x_{1},y_{1})
                                 \leq \cdots \leq (x_{m},y_{m}) = \ho}
         \cdots
           G(\Psi([(x_{i-1},y_{i-1}),(x_{i},y_{i})]))   
         \cdots
           t_{i}^{\rho((x_{i-1},y_{i-1}),(x_{i},y_{i}))}
         \cdots                             \\
  & = &
       \lim_{m \longrightarrow \infty}
\left(
       \sum_{\hz = x_{0} \leq x_{1} \leq \cdots \leq x_{m} = \ho}
         \cdots
           G(\Psi([x_{i-1},x_{i}]))   
         \cdots
           t_{i}^{\rho(x_{i-1},x_{i})}
         \cdots
\right)  \\
  &   &
\hspace*{20 mm}
\cdot
\left(
       \sum_{\hz = y_{0} \leq y_{1} \leq \cdots \leq y_{m} = \ho}
         \cdots
           G(\Psi([y_{i-1},y_{i}]))   
         \cdots
           t_{i}^{\rho(y_{i-1},y_{i})}
         \cdots
\right) \\
  & = &
      \widetilde{g}(F(P)) \cdot \widetilde{g}(F(Q)) ,
\end{eqnarray*}
where we only write the generic factor in each term.
\end{proof}

Combining
equation~(\ref{equation_Delta_widetilde_g})
in Proposition~\ref{proposition_g_widetilde_g}
with Theorem~\ref{theorem_widetilde_M},
we obtain:
\begin{theorem}
Let $G$ be a character on $\kab$.
Then the associated function $\widetilde{g}$ is a Hopf algebra endomorphism
on the quasisymmetric functions $\QSym$.
\end{theorem}

This theorem is a special case of a more general theorem due to
Aguiar, Bergeron and Sottile~\cite{Aguiar_Bergeron_Sottile}.
They proved that in 
the category of combinatorial
Hopf algebras
the quasisymmetric functions $\QSym$
is a terminal object.
A combinatorial Hopf algebra is
a Hopf algebra $H$ together with a character $G$.
Their results then states that
given a combinatorial Hopf algebra $H$ with character
$G$, there exists a Hopf algebra homomorphism
$\psi : H \longrightarrow \QSym$ such that
$G = \zeta \circ \psi$, where $\zeta$ is the character
on $\QSym$ defined by $\zeta(f) = A(\gamma^{-1}(f))$
and $\zeta(1_{\QSym}) = 1$.

\section{Type $B$ quasisymmetric functions}
\label{section_quasisymmetric_type_B}
\setcounter{equation}{0}

We now turn our attention to the chain map of the first kind.
In this section we will assume the underlying map $G$
is multiplicative with respect to the Cartesian product of posets.
The purpose of this section is to prove the 
chain map of the first kind $g$ is an algebra
map under the product $N$, and moreover,
to prove  $g$ is a comodule
map.
\begin{theorem}
For all $\ab$-polynomials $u$ and $v$, we have
\begin{eqnarray*}
      g(N(u,v))
  & = &
      N(g(u),
        g(v))      .
\end{eqnarray*}
\label{theorem_g_N}
\end{theorem}

By observing  
$N(\av \cdot u,  \av \cdot v) = \av \cdot M(u,v)$
(see Proposition~7.8 in~\cite{Ehrenborg_Fox}),
we have the corollary:
\begin{corollary}
For all $\ab$-polynomials $u$ and $v$, we have
\begin{eqnarray*}
      g(\av \cdot M(u,v))
  & = &
      N(g(\av \cdot u),
        g(\av \cdot v))      .
\end{eqnarray*}
\end{corollary}
This corollary implies 
Theorem~\ref{theorem_widetilde_M}
by applying the $H$ map.

In order to prove 
Theorem~\ref{theorem_g_N},
we introduce the type $B$ quasisymmetric functions
due to Chow~\cite{Chow}.
Let $\BQSym$ denote the algebra
$\kk[s] \tensor \QSym$. We view
$\BQSym$ as a subalgebra of
$\kk[s, t_1, t_2, \ldots] \cong \kk[s] \tensor \kk[t_1, t_2, \ldots]$.

Define the {\em type $B$ quasisymmetric function} of a poset $P$ by
\begin{eqnarray*}
     F_{B}(P)
  & = &
     \sum_{\hz < x \leq \ho}
         s^{\rho(x)-1} \cdot F([x,\ho]) \\
  & = &
     \lim_{m \longrightarrow \infty}
        \sum_{\hz < x_{0} \leq x_{1} \leq \cdots \leq x_{m} = \ho}
             s^{\rho(\hz,x_{0})-1}
           \cdot
             t_{1}^{\rho(x_{0},x_{1})}
           \cdot
             t_{2}^{\rho(x_{1},x_{2})}
           \cdots
             t_{m}^{\rho(x_{m-1},x_{m})}    .
\end{eqnarray*}

\begin{theorem}
For two graded posets $P$ and $Q$, we have 
$$        F_{B}(P \diamond Q)
      = 
          F_{B}(P) \cdot F_{B}(Q)   .  $$
\label{theorem_F_B_P_diamond_Q}
\end{theorem}
\begin{proof}
Applying the definition of $F_{B}$ to
the diamond product $P \diamond Q$ yields
\begin{eqnarray*}
     F_{B}(P \diamond Q)
  & = &
     \sum_{\hz < (x,y) \leq \ho_{P \diamond Q}}
         s^{\rho_{P \diamond Q}((x,y))-1}
      \cdot
         F([(x,y),\ho_{P \diamond Q}]) \\
  & = &
\left(
     \sum_{\hz < x \leq \ho_{P}}
         s^{\rho_{P}(x)-1}
      \cdot
         F([x,\ho_{P}])
\right)
  \cdot
\left(
     \sum_{\hz < y \leq \ho_{Q}}
         s^{\rho_{Q}(y)-1}
      \cdot
         F([y,\ho_{Q}])
\right) \\
  & = &
     F_{B}(P) \cdot F_{B}(Q)  .
\end{eqnarray*}
Here we are using 
$\rho_{P \diamond Q}((x,y))
  =
\rho_{P}(x)
  +
\rho_{Q}(y)
  - 1$ and
that
the quasisymmetric function $F$ is multiplicative
on posets. 
\end{proof}

Let $\gamma_{B}$ be the isomorphism between
$\kab$ and $\BQSym$ defined by
$$ \gamma_{B}\left(
        (\av-\bv)^{p}
      \cdot
        \bv
      \cdot
        (\av-\bv)^{p_{1}-1}
      \cdot
        \bv
      \cdots
        \bv
      \cdot
        (\av-\bv)^{p_{k}-1} \right)
 =
        s^{p}
      \cdot
        M_{(p_{1}, \ldots, p_{k})}     , $$
where $p \geq 0$ and $p_{1}, \ldots, p_{k} \geq 1$,
that is,
$\gamma_{B}(\Psi(P)) = F_{B}(P)$.
Define the linear map $g$ on $\BQSym$ by
$g(f) = \gamma_{B}(g(\gamma_{B}^{-1}(f)))$.
Hence Theorem~\ref{theorem_F_B_P_diamond_Q}
states 
\begin{equation}
 \gamma_{B}(N(u,v)) = \gamma_{B}(u) \cdot \gamma_{B}(v) .
\end{equation}

\begin{proposition}
The chain map of the first kind $g$
has the following form when applied to
the $\ab$-index, respectively, the type $B$ quasisymmetric function
of a poset $P$:
\begin{eqnarray}
g(\Psi(P))
  & = &
       \sum_{c}
           G(\Psi([x_{0},x_{1}]))
         \cdots
           G(\Psi([x_{k-1},x_{k}]))
         \cdot
           \wt(c)                           ,   \\
g(F_{B}(P))
  & = &
       \sum_{c}
           G(\Psi([x_{0},x_{1}]))
         \cdots
           G(\Psi([x_{k-1},x_{k}]))   
         \cdot
           s^{\rho(\hz,x_{0})-1}
         \cdot
           M_{(\rho(x_{0},x_{1}), \ldots, 
               \rho(x_{k-1},x_{k}))}      ,        \\
g(F_{B}(P))
  & = &
       \sum_{\hz < x \leq \ho}
           s^{\rho(x)-1}
         \cdot
           \widetilde{g}(F([x,\ho]))  , 
\label{equation_g_F_B}  \\
g(F_{B}(P))
  & = &
       \lim_{j \longrightarrow \infty}
       \sum_{m}
           G(\Psi([x_{0},x_{1}]))
         \cdots
           G(\Psi([x_{j-1},x_{j}]))     \nonumber \\
  &   &
\hspace*{40 mm}
         \cdot
           s^{\rho(x_{0})-1}
         \cdot
           t_{1}^{\rho(x_{0},x_{1})}
         \cdots
           t_{j}^{\rho(x_{j-1},x_{j})}    ,
\end{eqnarray}
where the first two sums
are over all chains
$c = \{\hz = x_0 < x_1 < \cdots < x_k = \ho\}$
in the poset $P$ and the fourth sum is over all
chains satisfying
$m = \{\hz < x_{0} \leq x_{1} \leq \cdots \leq x_{j} = \ho\}$
in $P$.
\label{proposition_g_Psi_F_B}
\end{proposition}
\begin{proof}
By the definition of $g$ and $\widetilde{g}$, we
have
\begin{eqnarray*}
g(u)
  & = &
\kappa(u)
  +
\sum_{u} \kappa(u_{(1)}) \cdot \bv \cdot \widetilde{g}(u_{(2)}) . 
\end{eqnarray*}
Apply this identity to the $\ab$-index of a poset $P$
and
use equation~(\ref{equation_widetilde_g_Psi})
to expand the factor $\widetilde{g}(\Psi([x,\ho]))$.
We then obtain
\begin{eqnarray*}
g(\Psi(P))
  & = &
(\av-\bv)^{\rho(P)-1}
  +
\sum_{\hz < x < \ho}
   (\av-\bv)^{\rho(x)-1}
        \cdot \bv \cdot
   \widetilde{g}(\Psi([x,\ho])) \\
  & = &
(\av-\bv)^{\rho(P)-1}
  +
\sum_{\hz < x < \ho}
       \sum_{k \geq 1}
       \sum_{x = x_{0} < x_{1} < \cdots < x_{k} = \ho}
           G(\Psi([x_{0},x_{1}]))
         \cdots
           G(\Psi([x_{k-1},x_{k}])) \\
  &   &
\hspace*{35 mm}
   \cdot (\av-\bv)^{\rho(x)-1}
        \cdot \bv
         \cdot
           (\av-\bv)^{\rho(x_{0},x_{1})-1}
         \cdot
           \bv
         \cdots
           \bv
         \cdot
           (\av-\bv)^{\rho(x_{k-1},x_{k})-1}         ,
\end{eqnarray*}
which is equivalent to the first identity of the proposition.

To prove the second identity,
apply the isomorphism $\gamma_{B}$ to the first identity.
The  third identity
follows from the second and
equation~(\ref{equation_widetilde_g_F_P_chain}).
Similarly, the fourth identity follows by
the third and equation~(\ref{equation_widetilde_g_F_P})
\end{proof}

As a remark,
we did not use the fact that
$G$ is a character on $\kab$
in the proof of Proposition~\ref{proposition_g_Psi_F_B}.

\begin{theorem}
The linear map $g$ is an algebra homomorphism
on the type $B$ quasisymmetric functions $\BQSym$.
\label{theorem_algebra}
\end{theorem}
\begin{proof}
By equation~(\ref{equation_g_F_B}) we have
\begin{eqnarray*}
     g(F_{B}(P \diamond Q))
  & = &
     \sum_{\hz < (x,y) \leq \ho_{P \diamond Q}}
         s^{\rho_{P \diamond Q}((x,y))-1}
      \cdot
         \widetilde{g}(F([(x,y),\ho_{P \diamond Q}])) \\
  & = &
\left(
     \sum_{\hz < x \leq \ho_{P}}
         s^{\rho_{P}(x)-1}
      \cdot
         \widetilde{g}(F([x,\ho_{P}]))
\right)
  \cdot
\left(
     \sum_{\hz < y \leq \ho_{Q}}
         s^{\rho_{Q}(y)-1}
      \cdot
         \widetilde{g}(F([y,\ho_{Q}]))
\right) \\
  & = &
     g(F_{B}(P)) \cdot g(F_{B}(Q))  .
\end{eqnarray*}
Since the type $B$ quasisymmetric function of posets
span the space $\BQSym$, the result follows.
\end{proof}

We remark that Theorems~\ref{theorem_g_N}
and~\ref{theorem_algebra} are equivalent via
the isomorphism $\gamma_{B}$.

Define the coproduct
$\Delta^{\BQSym} : \BQSym \longrightarrow \BQSym \tensor \QSym$
by
$$   \Delta^{\BQSym}(s^{p} \cdot f)
         =
      \sum_{f}^{\QSym} s^{p} \cdot f_{(1)} \tensor f_{(2)} , $$
where $f \in \QSym$. We then have the following result.

\begin{theorem}
For a graded poset $P$, we have
$$   \Delta^{\BQSym}(F_{B}(P))
   =
     \sum_{\hz < x \leq \ho}
           F_{B}([\hz,x])
         \tensor
           F([x,\ho])       .   $$
\label{theorem_Delta_QSym_B_F_B_P}
\end{theorem}
\begin{proof_qed}
By applying the definition of $F_{B}$, we obtain
\begin{eqnarray*}
\hspace*{15 mm}
\Delta^{\BQSym}(F_{B}(P))
 & = &
\Delta^{\BQSym}
\left(
  \sum_{\hz < y \leq \ho}
           s^{\rho(y)-1}
        \cdot
           F([y,\ho]) 
\right)            \\
 & = &
  \sum_{\hz < y \leq \ho}
  \sum_{y \leq x \leq \ho}
           s^{\rho(y)-1}
        \cdot
           F([y,x]) 
        \tensor
           F([x,\ho])     \\
 & = &
  \sum_{\hz < x \leq \ho}
\left(
  \sum_{0 < y \leq x}
           s^{\rho(y)-1}
        \cdot
           F([y,x]) 
\right)
        \tensor
           F([x,\ho])       \\
 & = &
  \sum_{\hz < x \leq \ho}
           F_{B}([\hz,x]) 
        \tensor
           F([x,\ho])     .
\hspace*{30 mm}
\hspace*{15 mm}
\qed
\end{eqnarray*}
\end{proof_qed}

\begin{theorem}
The linear map $g$ is a comodule endomorphism
on $\BQSym$, that is,
$$       \Delta^{\BQSym} \circ g 
      =
           (g \tensor \widetilde{g}) \circ \Delta^{\BQSym} . $$
\label{theorem_g_comodule_endomorphism}
\end{theorem}
\begin{proof}
By applying equation~(\ref{equation_g_F_B}) twice
and
Theorem~\ref{theorem_Delta_QSym_B_F_B_P},
we obtain
\begin{eqnarray*}
\Delta^{\BQSym}(g(F_{B}(P)))
  & = &
\Delta^{\BQSym}
\left(
       \sum_{\hz < x \leq \ho}
           s^{\rho(x)-1}
         \cdot
           \widetilde{g}(F([x,\ho]))
\right) \\
  & = &
       \sum_{\hz < x \leq \ho}
       \sum_{x \leq y \leq \ho}
           s^{\rho(x)-1}
         \cdot
           \widetilde{g}(F([x,y]))
         \tensor
           \widetilde{g}(F([y,\ho]))     \\
  & = &
       \sum_{\hz < y \leq \ho}
\left(
       \sum_{\hz < x \leq y}
           s^{\rho(x)-1}
         \cdot
           \widetilde{g}(F([x,y]))
\right)
         \tensor
           \widetilde{g}(F([y,\ho]))     \\
  & = &
       \sum_{\hz < y \leq \ho}
           g(F_{B}([\hz,y]))
         \tensor
           \widetilde{g}(F([y,\ho]))            \\
  & = &
       (g \tensor \widetilde{g}) 
\left(
       \sum_{\hz < y \leq \ho}
           F_{B}([\hz,y])
         \tensor
           F([y,\ho])
\right)                                            \\
  & = &
       (g \tensor \widetilde{g}) 
\left(
   \Delta^{\BQSym}(F_{B}(P))
\right)                                    .
\end{eqnarray*}
The result follows since 
the type $B$ quasisymmetric functions $F_{B}(P)$ span
$\BQSym$ as $P$ ranges over all posets.
\end{proof}

\section{Concluding remarks}
\setcounter{equation}{0}

Recall the following theorem of Hetyei~\cite[Theorem 1.10]{Hetyei_matrices}.
\begin{theorem}
If $P$ is the face poset of a spherical 
$CW$-complex
then
the Tchebyshev transform $T(P)$ is also
the face poset of a spherical
$CW$-complex.
\end{theorem}
A natural conjecture to make is
the following.
\begin{conjecture}
If $P$ is a spherical and shellable poset then
the Tchebyshev transform
$T(P)$ is also a spherical and shellable poset.
\end{conjecture}
See the related conjecture~\cite[Conjecture A.2]{Hetyei_matrices}.

A Gorenstein* lattice is an Eulerian
lattice which is Cohen-Macaulay.
See~\cite{Stanley_green_book} for terminology.
A very difficult question 
to settle is the Stanley's 
Gorenstein* conjecture~\cite{Stanley_Eulerian}:
among all Gorenstein* lattices
of rank $d+1$,
the $\cd$-index is minimized
on the 
simplex of dimension $d$.
This conjecture has been settled in the
special case of face
lattices polytopes
by Billera and Ehrenborg~\cite{Billera_Ehrenborg}.

Very little is known about Gorenstein* lattices.
Hetyei has proposed that the Tchebyshev transform
will give a host of new examples of Gorenstein* posets.
He has conjectured the following.
\begin{conjecture}
If $P$ is Gorenstein* then
the Tchebyshev transform $T(P)$ is Gorenstein*.
\end{conjecture}

One could also ask 
to study the behavior of
a Cohen-Macaulay poset $P$
under
the Tchebyshev transform $T$.
For example, can the system of parameters
and a basis for $T(P)$ be determined from the
original system of parameters and basis of
$P$.

Since the classical quasisymmetric functions
correspond to the symmetric group, that is,
the Weyl group of type $A$, the following
two questions are natural.
Are there analogues
of the quasisymmetric functions for other Weyl groups
other than type $A$ and $B$?
Similarly, are there analogues of the two maps
$F$ and $F_{B}$ on posets for the other Weyl groups?

For instance,
one can introduce one more extension of
the quasisymmetric function of a poset,
namely, define
$$
     F'(P)
  =
     \sum_{\hz < x \leq y < \ho}
         s^{\rho(\hz,x)-1} \cdot F([x,y]) \cdot u^{\rho(y,\ho)-1} .
$$
This poset invariant is
multiplicative with respect to the product
$(P - \{\hz,\ho\}) \times (Q - \{\hz,\ho\}) \cup \{\hz,\ho\}$
and has a bi-comodule structure.
Also, it behaves nicely with the map defined by
$$
     g'(u)
  =
       \sum_{k \geq 2}
       \sum_{u}
           \kappa(u_{(1)})
         \cdot
           \bv
         \cdot
           \widehat{g}(u_{(2)})
         \cdot
           \bv
         \cdots
           \bv
         \cdot
           \widehat{g}(u_{(k-1)})
         \cdot
           \bv
         \cdot
           \kappa(u_{(k)})   .  
$$
The essential question to answer is if these maps
naturally appear in geometry or combinatorics.
Also, is this a $\widetilde{B}$ analogue of the quasisymmetric
function of a poset?

Is there a notion of a type $B$ combinatorial Hopf algebra?
Moreover, is there an Aguiar, Bergeron and Sottile type theorem,
that is, that the pair $(\QSym,\BQSym)$ is the terminal object
in this category? These two question also extend to the other Weyl groups.

Another result due to Aguiar, Bergeron and
Sottile~\cite{Aguiar_Bergeron_Sottile}
is that every character $G$ of a Hopf algebra
factors into an even character $G_{+}$ and an odd character $G_{-}$.
In a recent preprint
Aguiar and Hsiao~\cite{Aguiar_Hsiao}
described this factorization explicitly
for the character $\zeta$.
The character $\zeta$ is the character
underlying 
Example~\ref{examples_four} (i).
Are there similar explicit factorizations
into even and odd characters
for Examples~\ref{examples_four} (ii) through (iv)?

We end with three open questions about the chain
maps $g$ and $\tilde{g}$.
Find other examples
of poset transformations so that  the resulting 
linear transformation on the $\ab$-index
has the form of $g$ or $\widetilde{g}$.
Find the general theorem
which determines the spectrum of the maps $g$ and 
$\widetilde{g}$.  
Alternatively, find subclasses of multiplicative maps  where this
is possible.
Recall for the Tchebyshev transform of the second
kind $\UTch$ we were able to do this.

\section*{Acknowledgements}

We graciously thank G\`abor Hetyei for inspiring
us to study the Tchebyshev transform
and suggesting research directions.
We thank Ira Gessel for directing us
to Chak-On Chow's work on
the type $B$ quasisymmetric functions.
The first author was partially supported by
National Science Foundation grant 0200624
and by 
a
University of Kentucky College of Arts~\& Sciences 
Faculty Research Fellowship.
The second author was partially supported by
a 
University of Kentucky College of Arts~\& Sciences Research Grant.
Both authors thank the Institute for Advanced Study/Park City
Mathematics Institute for providing a stimulating work environment.

\newcommand{\journal}[6]{{\sc #1,} #2, {\it #3} {\bf #4} (#5), #6.}
\newcommand{\book}[4]{{\sc #1,} ``#2,'' #3, #4.}
\newcommand{\bookf}[5]{{\sc #1,} ``#2,'' #3, #4, #5.}
\newcommand{\thesis}[4]{{\sc #1,} ``#2,'' Doctoral dissertation, #3, #4.}
\newcommand{\springer}[4]{{\sc #1,} ``#2,'' Lecture Notes in Math.,
                          Vol.\ #3, Springer-Verlag, Berlin, #4.}
\newcommand{\preprint}[3]{{\sc #1,} #2, preprint #3.}
\newcommand{\preparation}[2]{{\sc #1,} #2, in preparation.}
\newcommand{\appear}[3]{{\sc #1,} #2, to appear in {\it #3}}
\newcommand{\submitted}[4]{{\sc #1,} #2, submitted to {\it #3}, #4.}
\newcommand{\JCTA}{J.\ Combin.\ Theory Ser.\ A}
\newcommand{\AdvancesinMathematics}{Adv.\ Math.}
\newcommand{\JournalofAlgebraicCombinatorics}{J.\ Algebraic Combin.}

\newcommand{\communication}[1]{{\sc #1,} personal communication.}

{\em
\noindent
R.\ Ehrenborg,
Department of Mathematics,
University of Kentucky,
Lexington, KY 40506, \newline
{\tt jrge@ms.uky.edu}
}

{\em
\noindent
M.\ Readdy,
Department of Mathematics,
University of Kentucky,
Lexington, KY 40506, \newline
{\tt readdy@ms.uky.edu}
}

\end{document}